\documentclass[11pt]{article}
\usepackage[top=1in,bottom=1in,left=1in,right=1in]{geometry}

\usepackage{amsthm}
\usepackage{graphicx}
\usepackage{epstopdf}
\usepackage{algorithmic}
\usepackage{amsmath}
\usepackage{amssymb}
\usepackage{hyperref}
\usepackage[T1]{fontenc}
\usepackage{soul,xcolor}
\usepackage{cases}
\usepackage{booktabs}
\usepackage{float}
\usepackage{bm}
\usepackage{setspace}
\usepackage{makecell}
\usepackage{mathtools}
\usepackage{enumerate,enumitem}
\usepackage{multirow}
\usepackage{appendix}
\usepackage{natbib}
\usepackage{ulem}
\usepackage{xcolor}

\DeclareMathOperator{\GRAD}{GRAD}

\newcommand*{\dif}{\mathop{}\!\mathrm{d}}
\newtheorem{theorem}{Theorem}[section]

\newtheorem{lemma}[theorem]{Lemma}

\newtheorem{definition}{Definition}[section]

\title{MP-FVM: Enhancing Finite Volume Method for Water Infiltration Modeling in Unsaturated Soils via Message-passing Encoder-decoder Network}

\author{
 Zeyuan Song, Zheyu Jiang \\
  School of Chemical Engineering\\
  Oklahoma State University\\
  Stillwater, OK 74074 \\
  \texttt{\{taekwon.song,zjiang\}@okstate.edu} 
}

\begin{document}
\date{}
\maketitle

\begin{abstract}
The spatiotemporal water flow dynamics in unsaturated soils can generally be modeled by the Richards equation. To overcome the computational challenges associated with solving this highly nonlinear partial differential equation (PDE), we present a novel solution algorithm, which we name as the MP-FVM (Message Passing-Finite Volume Method), to holistically integrate adaptive fixed-point iteration scheme, encoder-decoder neural network architecture, Sobolev training, and message passing mechanism in a finite volume discretization framework. We thoroughly discuss the need and benefits of introducing these components to achieve synergistic improvements in accuracy and stability of the solution. We also show that our MP-FVM algorithm can accurately solve the mixed-form $n$-dimensional Richards equation with guaranteed convergence under reasonable assumptions. Through several illustrative examples, we demonstrate that our MP-FVM algorithm not only achieves superior accuracy, but also better preserves the underlying physical laws and mass conservation of the Richards equation compared to state-of-the-art solution algorithms and the commercial HYDRUS solver.


\end{abstract}

\section{Introduction}

The spatiotemporal dynamics of root zone (e.g., top 1 m of soil) soil moisture from precipitation and surface soil moisture information can generally be modeled by the Richards equation \citep{richards}, which captures irrigation, precipitation, evapotranspiration, runoff, and drainage dynamics in soil:
\begin{equation}\label{eqn_richards}
    \begin{aligned}
        &\partial_t\theta(\psi)+\nabla \cdot \textbf{q} = -S(\psi),\\
        &\textbf{q} = -K(\theta(\psi))\nabla (\psi+z).
    \end{aligned}
\end{equation}

Here, $\psi$ stands for pressure head (in, e.g., m), $\textbf{q}$ represents the water flux (in, e.g., $\text{m}^3/\text{m}^2\cdot \text{s}$), $S$ is the sink term associated with root water uptake (in, e.g., $\text{s}^{-1}$), $\theta$ denotes the soil moisture content (in, e.g., $\text{m}^3/\text{m}^3$), $K$ is unsaturated hydraulic water conductivity (in, e.g., $\text{m}/\text{s}$), $t \in [0,T]$ denotes the time (in, e.g., $\text{s}$), and $z$ corresponds to the vertical depth (in, e.g., $\text{m}$). The Richards equation is a nonlinear convection-diffusion equation \citep{caputo_front_2008}, in which the convection term is due to gravity, and the diffusive term comes from Darcy's law \citep{smith_infiltration_2002}. For unsaturated flow, both $\theta$ and $K$ are highly nonlinear functions of pressure head $\psi$ and soil properties, making Equation \eqref{eqn_richards} challenging to solve numerically \citep{escape_25_inverse}. Specifically, $\theta(\psi)$ and $K(\psi)$ (or $K(\theta)$, depending on the model) are commonly referred to as the water retention curve (WRC) and hydraulic conductivity function (HCF), respectively. Several of the most widely used empirical models for WRC and HCF are summarized in Table \ref{table_hcfwrc}.

\begin{table}[ht!]
    \centering
    \begin{tabular}{l l l}
    \toprule
    Model & HCF ($K(\psi)$ or $K(\theta)$) & WRC ($\theta(\psi)$)\\
    \midrule
    \citet{haverkamp1977comparison} & $K_s\frac{A}{A+|\psi|^{\gamma}}$ & $\theta_r+\frac{\alpha(\theta_s-\theta_r)}{\alpha+|\psi|^\beta}$ \\
    \citet{Mualem1976,vanGenuchten1980} &
    $K_s\sqrt{\frac{\theta-\theta_r}{\theta_s-\theta_r}} \left\{1-\left[1-\left(\frac{\theta-\theta_r}{\theta_s-\theta_r}\right)^{\frac{l}{l-1}}\right]^{\frac{l-1}{l}}\right\}^2$
    & $\theta_r+\frac{\theta_s-\theta_r}{\big[1+(\alpha|\psi|)^n\big]^{\frac{n-1}{n}}}$ \\
    \citet{gardner1958some} & $K_s e^{\alpha \psi}$ & $\theta_r + (\theta_s-\theta_r)e^{\alpha \psi}$\\
    \bottomrule
    \end{tabular}
    \caption{Some of the widely used HCF and WRC models. In these models, $A$, $\gamma$, $\alpha$, $\beta$, $n$, $\theta_s$, and $\theta_r$ are soil-specific parameters and have been tabulated for major soil types.}
    \label{table_hcfwrc}
\end{table}

Due to the highly nonlinear nature of WRC and HCF, analytical solutions to the Richards equation do not exist in general \citep{farthing2017numerical}. Thus, the Richards equation is typically solved numerically in some discretized form \citep{songescape}. Consider the discretized version of Equation \eqref{eqn_richards}, whose control volume $V \subset \mathbb{R}^d$ ($d=1,2,3$) is discretized into $N$ small cells $V_1, \dots, V_N$. Using implicit Euler method on the time domain with a time step size of $\Delta t$, the discretized Richards equation at time step $m = 0,1, \dots, \lceil\frac{T}{\Delta t}\rceil-1$ can be expressed as:
\begin{equation}\label{eqn_discretizedRE}
    \left\{
    \begin{aligned}
        &\theta(\psi_i^{m+1})-\theta(\psi_i^{m})-\Delta t \nabla \cdot \Big[K\left(\theta(\psi_i^{m+1})\right)\nabla\left(\psi_i^{m+1}+z\right)\Big] + \Delta tS(\psi_i^{m+1})= 0,\\
        &\text{Dirichlet boundary condition: } \psi_j (\cdot)=0 \quad \text{for all } V_j \subset \partial V,\\
        &\text{Initial condition: }\psi(0,\cdot)=\psi_0(\cdot),
    \end{aligned}
    \right.
\end{equation}
where $\psi_i^{m}$ is the pressure head in cell $V_i$ and time step $m$, and $\psi_0(\cdot)$ denotes the initial condition at $t=0$. 

The performance of a numerical PDE solver depends theoretically on the well-posedness of the PDE \citep{sizikov2011well}, which is an essential property that certifies the accuracy and reliability of numerical solutions to the PDE. A PDE is said to be well-posed if its weak solution exists, is unique, and depends continuously on the problem's initial conditions \citep{sizikov2011well,evans2010partial}. Here, we consider an FVM discretization with a discrete space $Q_h\subset L^2(V)$ of piecewise constants, where $h$ denotes the maximum dimension of any cell in its mesh. With this, we define the space of piecewise constant functions on the set of meshes $\mathcal{T}_h = \{V_1, V_2, \ldots, V_N\}$ as $Q_h(V) = \{ v \in L^2(V) : v|_{V_i} \text{ is constant for all } V_i \in \mathcal{T}_h \}$. Then, we introduce the discrete gradient operator \citep{hyman1997natural}, $\text{GRAD}_h$, which maps a cell-based function in $Q_h$ to a face-based function that approximates the gradient. Note that $\psi_i^m$ in Equation \eqref{eqn_discretizedRE} denotes the pressure head in cell $V_i$ and time step $m$, which is the value of $\psi^m$ in the cell $V_i$. To study the pressure head solution in function space $Q_h(V)$, we focus on $\psi^m$ rather than  $\psi_i^m$. With this, the discrete solution for the FVM-discretized Richards equation can be defined as follows:
\begin{definition}\label{def_weaksolution}
 Given $\psi^{m}\in Q_h,$ if for any $v\in Q_h$,
\begin{equation}\label{eqn_definition}
    \left\langle\theta(\psi^{m+1})-\theta(\psi^{m}),v\right\rangle_{V} + \Delta t\Big\langle K(\theta(\psi^{m+1}))\GRAD_h(\psi^{m+1}+z), \GRAD_h (v) \Big\rangle_{\mathcal{E}_h} + \Delta t\left\langle S(\psi^{m+1}), v\right\rangle_{V}=0
\end{equation}
holds, where $\mathcal{E}_h$ denotes the set of all faces that make up the mesh $\mathcal{T}_h$, then $\psi^{m+1}$ is a discrete solution of the FVM-discretized Richards equation.
\end{definition}

Following Definition \ref{def_weaksolution}, for the discrete function space $Q_h$, an inner product over a cell $V_i$ is defined for piecewise constant functions $f, g \in Q_h$ as $\langle f,g \rangle_{V_i}\coloneqq \int_{V_i} fg \dif V$. In this case, by denoting $f_i$ and $g_i$ as the function values of $f$ and $g$ respectively on $V_i$ (i.e., $f_i=f|_{V_i}$ and $g_i=g|_{V_i}$, both of which are constants), we have $\int_{V_i} fg \dif V = f_ig_i \text{vol}(V_i)$. The global inner product over the entire domain $V$ is then $\langle f,g \rangle_{V} \coloneqq \sum_{i=1}^N \langle f,g \rangle_{V_i}$. We remark that the existence and uniqueness of the weak solution of the Richards equation have been rigorously established and carefully studied \citep{merz2010strong,misiats2013second,abdellatif2016priori}, setting up the theoretical foundation for developing an efficient solution algorithm to solve the discretized Richards equation numerically.

\section{Literature Review}

Among existing solution algorithms for the Richards equation, methods based on finite difference and finite element discretizations \citep{day1956numerical,celia1990general,chavez2024interface, haghighat2023finite} are the most well studied and implemented \citep{richardsreview}. However, these methods often face challenges when handling large-scale problems and suffer from instability issues such as oscillations \citep{oscillation}. Recently, \citet{ireson2023simple} used the method of lines to convert the 1-D Richards equation into an ordinary differential equation (ODE), which was then solved by finite difference method. Similarly, the process converting 1-D Richards equation into an ODE can also be achieved by implementing generalized Boltzmann transform \citep{zhou2013richards}. Despite these advancements, finite difference- and finite element-based methods generally require high mesh resolution to satisfy the local equilibrium condition \citep{or2015,roth2008,Vogel2008}. Furthermore, they tend to fail to preserve global mass conservation \citep{massconservative} and other important underlying physical relations among soil moisture, pressure head, and water flux.

Meanwhile, finite volume discretization method (FVM) has the potential to achieve high solution accuracy and preserve mass conservation when solving the Richards equation \citep{eymard1999finite}. Some of the notable works include \citet{lai2015mass} who obtained a family of mass-conservative finite volume predictor-corrector solutions for the 1-D Richards equation, \citet{misiats2013second} who proposed a second-order accurate monotone FVM for 1-D Richards equation, as well as others \citep{bassetto2022several, caviedes2013verification, manzini2004mass, su2022numerical}. However, like finite difference and finite element methods, conventional FVM typically converts the discretized Richards equation into a large, stiff matrix equation, which can be challenging to solve.

Instead of following the standard practice of converting the discretized equation into a matrix equation, a fixed-point iteration scheme solves the discretized Richards equations iteratively. Some of the notable fixed-point iteration schemes include \citet{bergamaschi1999mixed} and \citet{zeidler2013nonlinear}. In standard fixed-point iteration schemes, a static parameter $\tau$ is used for all time steps and discretized cells. However, choosing an appropriate $\tau$ value is not straightforward, as soil moisture content and pressure head can exhibit strong spatiotemporal variations. For instance, \citet{zhu2019improved} reported that, when a large static parameter $\tau$ is chosen, the discretized Richards equation could become ill-posed. To address this drawback, \citet{amrein2019adaptive} proposed a fully adaptive fixed-point iteration scheme based on the Galerkin method to solve 2-D semilinear elliptic equations by finding the optimal mesh refinements for fixed $\tau$ at each iteration. Nevertheless, some of the existing adaptive fixed-point iteration schemes have been reported to suffer from numerical oscillations \citep{casulli2010nested}. So far, no fully adaptive fixed-point iteration scheme has been established for solving the 3-D Richards equation.

With the breakthroughs in artificial intelligence and machine/deep learning, a new avenue for solving PDEs is to directly incorporate physical knowledge and constraints derived from the PDE into a neural network. One of the popular frameworks is the Physics-Informed Neural Network (PINN) \citep{raissi2019physics,raissi2017physics,chen2023modeling,lan2024reconstructing,ng2025novel}, where the PDE itself is embedded in the loss function as a regularization term. However, this often results in high computational costs and training instability. On the other hand, hybrid methods that integrate machine/deep learning techniques with discretization-based numerical methods have emerged as a promising approach to enhance the accuracy and stability of numerical algorithms \citep{bar2019learning}. In particular, the encoder-decoder network architecture has shown great potential in solving PDEs \citep{pichi2024graph,lu2020}. For example, \citet{ranade2021discretizationnet} employed a generative CNN-based encoder-decoder model with PDE variables as both input and output features. However, hybrid PDE solvers based on conventional encoder-decoder architecture face convergence and stability issues as the latent variables do not have physical meaning or obey conservation laws. To overcome this limitation, \citet{brandstetter} proposed a message passing mechanism, in which a new module called processor consisting of graph neural networks is placed between the encoder and the decoder. Meanwhile, \citet{ranade2021discretizationnet} proposed the use of iterative algorithm along with the encoder-decoder network, such that the latent variables can be viewed as the numerical solutions of PDE and convergence can be established without invoking the message passing mechanism. Nevertheless, to the best of our knowledge, the encoder-decoder network architecture and the message passing mechanism have not yet been integrated to solve the Richards equation.

\section{Motivation and Scope of Our Approach}

When solving the FVM-discretized Richards equation using static fixed-point iteration scheme, we observe that different choices of static parameter $\tau$ and total iteration count $S$ lead to slightly different solutions. This aligns with the observations made by \citet{zhu2019improved} that varying $\tau$ affects the condition number of the discretized Richards equation. To ensure our numerical scheme is well-posed and reaches convergence within $S$ iterations, here we propose an adaptive fixed-point iteration scheme, where $\tau$ can dynamically adjust itself with respect to space, time and iteration count. Furthermore, we propose to characterize the solution discrepancies in a data-driven approach \citep{songdata}. Specifically, we adopt the encoder-decoder network architecture to solve the Richards equation by extending the message passing mechanism to FVM-based adaptive fixed-point iteration scheme. This results in a novel numerical framework called the MP-FVM (Message Passing Finite Volume Method). In addition to these architectural advancements, we introduce several new strategies to further improve the numerical accuracy and computational efficiency of our MP-FVM algorithm, including:

\begin{itemize}[leftmargin=*]
    \item We present a generalized framework for any $d$-dimensional ($d=1,2,3$) Richards equation, and demonstrate that the new MP-FVM algorithm can be versatilely adopted to modeling different realistic scenarios (e.g., layered soil, actual precipitation). We provide rigorous theoretical justification of the convergence behavior of our MP-FVM algorithm.
    \item We introduce a ``coarse-to-fine'' approach to enhance the solution accuracy of our MP-FVM algorithm without requiring a large amount of high-accuracy, fine-mesh training data. We demonstrate that this coarse-to-fine approach maintains a good balance between computational efficiency and solution accuracy.
    \item We show that, by synergistically integrating our novel adaptive fixed-point iteration scheme, FVM, and the encoder-decoder network, our new MP-FVM algorithmic framework significantly enhances the ability of FVM discretization in preserving the underlying physical relationships and mass conservation associated with the Richards equation.
\end{itemize}

Overall, these techniques holistically improve convergence and reduce numerical oscillations compared to conventional FVM. Meanwhile, they also help preserve physical laws (e.g., global water balance more reliably compared to standard finite element and finite difference schemes. Furthermore, our MP-FVM algorithm achieves fine-scale accuracy using only coarse-grid training data, hence bypassing computationally expensive training. Finally, we remark that MP-FVM algorithm can leverage pre-trained models to reduce retraining time, enabling transfer learning across different boundary and/or initial conditions.
  
We organize the subsequent sections of the paper as follows. In Section \ref{section_FVM}, we derive the FVM-based adaptive fixed-point iteration formulation and prove its global convergence. Then, in Section \ref{sec_DRW}, we present the encoder-decoder network architecture as well as its integration with the message passing mechanism to form the MP-FVM algorithm. To compare the performance of our MP-FVM algorithm with state-of-the-art Richards equation solvers, we conduct comprehensive case studies and in-depth analyses of 1- through 3-D benchmark problems in Section \ref{sec_casestudies}. Finally, we summarize the results and discuss future directions in Section \ref{conclusion}.

\section{Adaptive fixed-point iteration scheme of Discretized Richards Equation}\label{section_FVM} 
In this section, we will formally introduce the adaptive fixed-point iteration scheme formulation of the FVM-discretized Richards equation. We will also derive sufficient conditions for parameter $\tau$ to ensure convergence. We will also analyze the convergence behavior of the resulting sequence of solutions $\{\psi_i^{m+1,s}\}_s$, where $s$ is the iteration count ($s=1,2,\ldots, S$).

\subsection{Adaptive fixed-point iteration scheme for the Richards Equation}

To discretize the Richards equation via FVM, we first integrate both sides of Equation \eqref{eqn_richards} over $V$:
\begin{equation}\label{eqn_integralRE}
    \int_V  \left[\partial_t\theta(\psi) + S(\psi) \right] \dif V=\int_V \nabla \cdot \big[K(\theta)\nabla (\psi+z)\big]\dif V.
\end{equation}

Next, we apply the divergence theorem to Equation \eqref{eqn_integralRE}, which converts the volume integral on the RHS into a surface integral:
 \begin{equation}\label{eqn_surfaceRE}
   \left[\partial_t\theta(\hat{\psi}) + S(\hat{\psi}) \right]_{\hat{\psi}\in V}\text{vol}(V)=\oint_{S_V} K(\theta)\nabla (\psi+z)\cdot\mathbf{n} \dif S_V,
\end{equation}
where $\text{vol}(V)$ is the volume of $V$, $S_V$ is the surface of $V$ and $\mathbf{n}$ is the outward pointing unit normal to the boundary $\partial V$. The common surface shared by cell $V_i$ and cell $V_j$ is denoted as $\omega_{i,j}$. With this, we can rewrite the operator $K(\cdot)\nabla (\cdot)$ and the outward pointing unit normal vector $\mathbf{n}$ on $\omega_{i,j}$ as $\big[K(\cdot)\nabla (\cdot)\big]_{\omega_{i,j}}$ and $\mathbf{n}_{\omega_{i,j}}$, respectively. After FVM discretization, we obtain the discretized version of Equation \eqref{eqn_surfaceRE} as:
\begin{equation} \label{eqn_fvm}
    \partial_t\theta_i\text{vol}(V_i) + S(\psi_i) \text{vol}(V_i) = \sum_{j\in \mathcal{N}_i} \big[K(\theta)\nabla (\psi+z)\big]_{\omega_{i,j}} \cdot\mathbf{n}_{\omega_{i,j}} A_{\omega_{i,j}} \qquad \forall i=1,\dots,N,
\end{equation}
where $\partial_t\theta_i$ refers to the time derivative $\partial_t\theta(\psi_i)$ in cell $V_i$, $\mathcal{N}_i$ denotes the index set of all the neighboring cells sharing a common surface with $V_i$, and $A_{\omega_{i,j}}$ is the area of surface $\omega_{i,j}$.

In static fixed-point iteration scheme, for each cell $V_i$ and at each time step $m+1$, one would add the term $\frac{1}{\tau}(\psi_i^{m+1,s+1} - \psi_i^{m+1,s})$ to either side of Equation \eqref{eqn_fvm}, so that the Richards equation can be solved in an iterative manner. The fixed-point pressure head solution of this iterative procedure is denoted as $\psi_i^m$. Since $\tau$ is a static constant, a trial-and-error procedure is typically required to obtain an appropriate $\tau$ value that avoids convergence issues. Not only is this search procedure tedious to implement, the solutions obtained are also less accurate most of the time as we will show in Section \ref{sec_1d}. Thus, inspired by previous works \citep{amrein2019adaptive,zhu2019improved}, we propose an adaptive fixed-point iteration scheme that replaces the static $\tau$ with $\tau_i^{m+1,s}$, which adjusts itself for each specific discretized cell, time step, and iteration count. We then introduce the term $\frac{1}{\tau_i^{m+1,s}}(\psi_i^{m+1,s+1}-\psi_i^{m+1,s})$ to the LHS of Equation \eqref{eqn_fvm}, which leads to:
\begin{equation}\label{eqn_lscheme1}
    \begin{aligned}
     \psi_i^{m+1,s+1} = & \psi_i^{m+1,s}+\tau_i^{m+1,s}\sum_{j\in \mathcal{N}_i} \big[K(\theta)\nabla (\psi+z)\big]_{\omega_{i,j}}^{m+1,s}\cdot\mathbf{n}_{\omega_{i,j}}A_{\omega_{i,j}} \\
     & - \tau_i^{m+1,s} \left[\partial_t\theta_i^{m+1,s} + S(\psi_i^{m+1,s})\right]\text{vol}(V_i),
     \end{aligned}
\end{equation}

By discretizing $\partial_t\theta_i^{m+1,s}$ using implicit Euler scheme as $\frac{\theta(\psi_i^{m+1,s})-\theta(\psi_i^m)}{\Delta t}$, we can obtain the adaptive fixed-point iteration scheme of the FVM-discretized Richards equation:
\begin{equation}\label{eqn_lscheme2}
    \begin{aligned}
         \psi_i^{m+1,s+1} &=  \psi_i^{m+1,s}+\tau_i^{m+1,s}\sum_{j\in \mathcal{N}_i} K_{\omega_{i,j}}^{m+1,s}\frac{(\psi+z)_{j}^{m+1,s}-(\psi+z)_{i}^{m+1,s}}{\text{dist}(V_j,V_i)} \mathbf{e} \cdot \mathbf{n}_{\omega_{i,j}}A_{\omega_{i,j}}\\
         & - \tau_i^{m+1,s} \left[\frac{\theta(\psi_i^{m+1,s})-\theta(\psi_i^m)}{\Delta t} + S(\psi_i^{m+1,s})\right]\text{vol}(V_i),
     \end{aligned}
\end{equation}
where $\mathbf{e}=(1,1,1)$ for the standard 3-D Cartesian coordinate system, and $\text{dist}(\cdot,\cdot)$ represents the Euclidean distance function.

\subsection{Choice of Adaptive Linearization Parameter}
In adaptive fixed-point iteration scheme, we observe that $\tau_i^{m+1,s}$ needs to be sufficiently small because otherwise, the RHS of Equation \eqref{eqn_lscheme2} could approach infinity, which affects the convergence of the scheme. To prevent $\frac{1}{\tau_i^{m+1,s}}$ from being too large, we impose a user-specified global upper bound $\tau_0$:
\begin{equation*}
    \tau_i^{m+1,s} \leq \tau_0.
\end{equation*}

In addition, the choice of $\tau_i^{m+1,s}$ can impact the accuracy of solutions. In other words, the term $\left|\frac{\psi_i^{m+1,s+1}-\psi_i^{m+1,s}}{\psi_i^{m+1,s}}\right|$ should be no greater than a prespecified tolerance $\rho$. Thus, we have:
\begin{equation*}
    \left|\frac{\psi_i^{m+1,s+1}-\psi_i^{m+1,s}}{\psi_i^{m+1,s}}\right|=\frac{\tau_i^{m+1,s}|g_i^{m+1,s}|}{|\psi_i^{m+1,s}|}\leq\rho \qquad \forall s=1,\dots,S,
\end{equation*}
where $S$ is the user-specified total number of iterations for convergence, $\rho$ should be no less than the overall tolerance of convergence $\epsilon$ (to be discussed in Section \ref{sec_convergence}), and: 
\begin{equation*}
    \begin{aligned}
        g_i^{m+1,s} = &\sum_{j\in \mathcal{N}_i} K_{\omega_{i,j}}^{m+1,s}\frac{(\psi+z)_{j}^{m+1,s}-(\psi+z)_{i}^{m+1,s}}{\text{dist}(V_j,V_i)} \mathbf{e} \cdot \mathbf{n}_{\omega_{i,j}}A_{\omega_{i,j}} \\
        & -\frac{\theta(\psi_i^{m+1,s})-\theta(\psi_i^m)}{\Delta t}\text{vol}(V_i) - S(\psi_i^{m+1,s}) \text{vol}(V_i).
    \end{aligned}
\end{equation*}

This implies that:
\begin{equation}\label{eqn_Lbound1}
    \tau_i^{m+1,s}\leq \frac{\rho|\psi_i^{m+1,s}|}{(1+\rho)|g_i^{m+1,s}|} \qquad \forall s=1,\dots,S,
\end{equation}
whose RHS can be explicitly determined from the results of the previous iteration. Note that, in actual implementation, we select $\tau_i^{m+1,s}$ based on:
\begin{equation}\label{eqn_Lbound2}
    \tau_i^{m+1,s} = \min\left\{\tau_0, \frac{\rho|\psi_i^{m+1,s}|}{(1+\rho)|g_i^{m+1,s}|}\right\} \qquad \forall s=1,\dots,S.
\end{equation}

Meanwhile, we can monitor the sensitivity of solutions obtained by our adaptive fixed-point iteration scheme and make sure that the solutions do not change drastically with respect to small perturbations. To achieve this, following \citet{zarba1988thesis} and \citet{zarba1988numerical}, we explicitly write down Equation \eqref{eqn_lscheme2} for all discretized cells in the form of a matrix equation:
\begin{equation}\label{eqn_matrix}
    \mathbf{A} \mathbf{x}^{m+1,s} = \mathbf{b},
\end{equation}
where the $i$th element of vector $\mathbf{x}^{m+1,s}$ is $x_i^{m+1,s} = \psi_i^{m+1,s+1} - \psi_i^{m+1,s}$, which corresponds to cell $V_i$. Here, it is worth mentioning that Equation \eqref{eqn_matrix} is not used for solving Equations \eqref{eqn_lscheme2} as it is an explicit numerical scheme. Rather, it is used for analyzing the properties of the scheme after $x_i^{m+1,s}$ solutions are obtained by solving Equation \eqref{eqn_lscheme2}. For example, to evaluate the choice of $\tau_i^{m+1,s}$, we can calculate the condition number of $\mathbf{A}$ based on the solutions obtained from the chosen $\tau_i^{m+1,s}$. If the condition number is larger than a user-specified threshold, we will update $\tau_0$ in Equation \eqref{eqn_Lbound2} so that the condition number drops below the threshold. For 1-D problems, \citet{zarba1988thesis} showed that $\mathbf{A}$ is a $N\times N$ asymmetric tridiagonal matrix. In this case, the condition number of $\mathbf{A}$ can be determined by calculating its eigenvalues. On the other hand, for 2-D and 3-D problems, $\mathbf{A}$ is a rectangular matrix, so that singular value decomposition will be used to determine its condition number.

\subsection{Convergence of Adaptive Fixed-Point Iteration Scheme}\label{sec_convergence}
We now study the convergence behavior of our adaptive fixed-point iteration scheme, which is formalized in Theorem \ref{theorem_convergence}. Recall that functions $\psi^{m}$ and $\psi^{m+1,s}$ are considered to study the convergence.

To show this, the idea is to leverage Definition \ref{def_weaksolution} and find $\psi^{m+1,s+1} \in Q_h(V)$ given $\psi^{m},\,\psi^{m+1,s} \in Q_h(V)$ such that:
\begin{equation}\label{eqn_discretizedRE2}
\begin{split}
    {\left\langle\theta(\psi^{m+1,s+1})-\theta(\psi^{m}),v\right\rangle}_{V} &+ \frac{\Delta t}{\tau^{m+1,s}}{\left\langle\psi^{m+1,s+1}-\psi^{m+1,s}, v\right\rangle}_{V} + {\left\langle S(\psi^{m+1,s+1}) , v \right\rangle}_{V} \\
    &=-\Delta t {\left\langle K\left(\theta(\psi^{m+1})\right)\text{GRAD}_h(\psi^{m+1,s+1}+z),\text{GRAD}_h(v)\right\rangle}_{\mathcal{E}_h}  
\end{split}
\end{equation}
holds for any $v \in Q_h(V)$. We remark that, unlike previous proofs (e.g., \citet{amrein2019adaptive}) that are based on several restrictive assumptions, our convergence proof follows a different approach that is intuitive and flexible, as it does not involve any additional assumptions other than the properties listed below.

\begin{theorem}\label{theorem_convergence}
    The sequence $\{\psi^{m+1,s}\}_s$ converges to a unique solution $\psi^{m+1} \in Q_h(V)$ for $m = 0,1,\ldots,\lceil\frac{T}{\Delta t}\rceil-1$.
\end{theorem}
\begin{proof}

First, we state two key properties used in the proof:
\begin{enumerate}[leftmargin=*,label= Observation \arabic*:]
\item The Cauchy-Schwarz inequality holds for the discrete $L^2$ inner product: for any $u,w \in Q_h$, we have $\langle u, w \rangle_V \le \|u\|_{L^2} \|w\|_{L^2}$.

\item $\dot{\theta}(\psi) = \frac{\dif \theta}{\dif \psi}|_{\psi^{m+1,s}} \ge c_0 > 0$, which is valid in most WRC models (see Table \ref{table_hcfwrc}). Similarly, $\dot{S}(\psi) = \frac{\dif S}{\dif \psi} |_{\psi^{m+1,s}} \geq 0 $ in the region between the start and optimal root water extraction.
\end{enumerate}

First, we subtract Equation \eqref{eqn_definition} from Equation \eqref{eqn_discretizedRE2} to obtain the error equation. Let $e^s \coloneqq \psi^{m+1,s} - \psi^{m+1}$, we have:
\begin{equation} \label{eqn_substractionRE_1}
\begin{split}
     & {\left\langle\theta(\psi^{m+1,s+1})-\theta(\psi^{m+1}),v\right\rangle}_{V} + \frac{\Delta t}{\tau^{m+1,s}}{\left\langle e^{s+1}-e^{s}, v \right\rangle}_{V}\\ 
     & + {\left\langle S(\psi^{m+1,s+1}) - S(\psi^{m+1}), v \right\rangle}_{V}
     = -\Delta t {\left\langle K(\cdot){\GRAD_h}(e^{s+1}),{\GRAD_h}(v)\right\rangle}_{{\mathcal{E}_h}}.  
\end{split}
\end{equation}

Let the test function $v=e^{s+1} = \psi^{m+1,s+1}-\psi^{m+1}$. This is a valid choice as $e^{s+1} \in Q_h$. By applying the mean value theorem to the $\theta$ and $S$ terms, and using Observation 2, we have:
\begin{equation}\label{eqn_nnn}
\begin{aligned}
    \left\langle\theta(\psi^{m+1,s+1})-\theta(\psi^{m+1}),e^{s+1}\right\rangle_{V} &= \left\langle \dot{\theta}(\xi_\theta) e^{s+1}, e^{s+1} \right\rangle_V \ge c_0 \|e^{s+1}\|_{L^2}^2, \\
    \left\langle S(\psi^{m+1,s+1})-S(\psi^{m+1}),e^{s+1}\right\rangle_{V} &= \left\langle \dot{S}(\xi_S) e^{s+1}, e^{s+1} \right\rangle_V \ge 0
\end{aligned}
\end{equation}
for some $\xi_\theta, \xi_S$ between  $\psi^{m+1,s+1}$ and $\psi^{m+1}$. The flux term on the RHS of Equation \eqref{eqn_substractionRE_1} is also non-negative:
\begin{equation}\label{eqn_nne}
  -\Delta t {\left\langle K(\cdot)\GRAD_h(e^{s+1}),{\GRAD_h}(e^{s+1})\right\rangle}_{{\mathcal{E}_h}} = -\Delta t \|{e^{s+1}}\|_{{h}}^2 \le 0,
\end{equation}
where $\|\cdot\|_h$ is the discrete energy semi-norm. Substituting Equation \eqref{eqn_nnn} and Equation \eqref{eqn_nne} into Equation \eqref{eqn_substractionRE_1} gives:
\begin{equation}\label{eqn_bb}
  c_0 \|e^{s+1}\|_{L^2}^2 + 0 + \frac{\Delta t}{\tau^{m+1,s}}{\left\langle e^{s+1}-e^{s}, e^{s+1} \right\rangle}_{V} \le 0  
\end{equation}
By applying Observation 1, Equation \eqref{eqn_bb} leads to:
\begin{equation}\label{eqn:bbb}
c_0 \|e^{s+1}\|_{L^2}^2 + \frac{\Delta t}{\tau^{m+1,s}} \left( \|e^{s+1}\|_{L^2}^2 - \langle e^s, e^{s+1} \rangle_V \right) \le 0. 
\end{equation}
Then, we have:
\begin{equation}\label{eqn:18}
\left( c_0 + \frac{\Delta t}{\tau^{m+1,s}} \right) \|e^{s+1}\|_{L^2}^2 \le \frac{\Delta t}{\tau^{m+1,s}} \langle e^s, e^{s+1} \rangle_V \le \frac{\Delta t}{\tau^{m+1,s}} \|e^s\|_{L^2} \|e^{s+1}\|_{L^2}. 
\end{equation}
If $e^{s+1}=0$, we complete the proof.
If $e^{s+1}\neq 0$, we can divide Equation 
\eqref{eqn:18} by $\|e^{s+1}\|_{L^2}$:
\begin{equation}
\left( c_0 + \frac{\Delta t}{\tau^{m+1,s}} \right) \|e^{s+1}\|_{L^2} \le \frac{\Delta t}{\tau^{m+1,s}} \|e^s\|_{L^2},     
\end{equation}
which yields the contraction:
\begin{equation}
    \|e^{s+1}\|_{L^2} \le \underbrace{\left(\frac{\frac{\Delta t}{\tau^{m+1,s}}}{c_0 + \frac{\Delta t}{\tau^{m+1,s}}}\right)}_{\eqqcolon \gamma_s} \|e^s\|_{L^2}.
\end{equation}
Since $c_0 > 0$, the contraction factor $\gamma_s$ is strictly less than 1. Therefore, the sequence is a contraction mapping on the discrete space $Q_h(V)$ equipped with the $L^2$ norm. By the Banach fixed-point theorem, the sequence $\{\psi^{m+1,s}\}$ converges to a unique solution $\psi^{m+1} \in Q_h(V)$. This completes the proof.
\end{proof}

\section{Message Passing Finite Volume Method (MP-FVM)} \label{sec_DRW}

\begin{figure}[ht!]
    \centering
    \includegraphics[width=\textwidth]{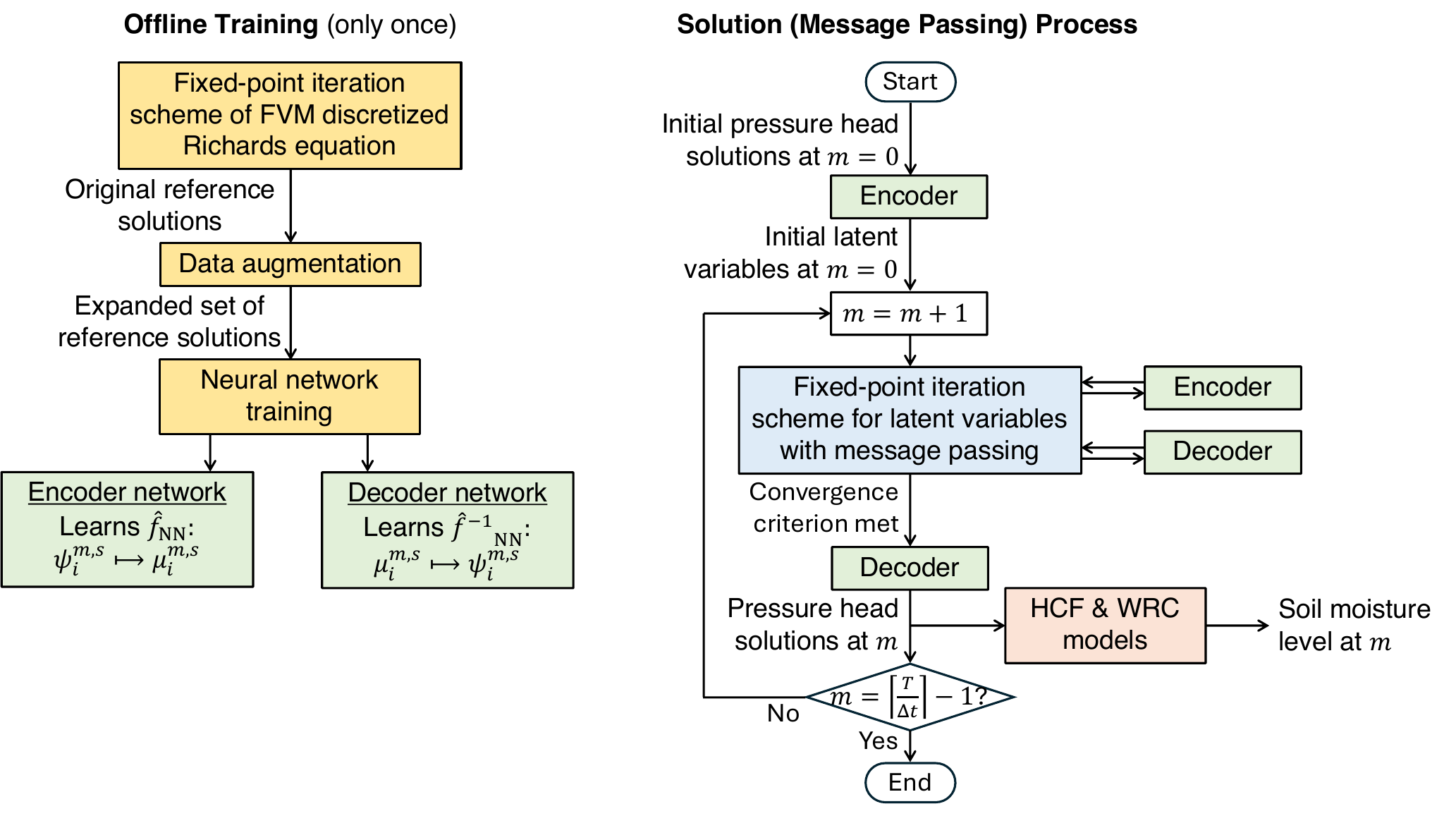}
    \vspace{-2em}
    \caption{Flowchart of our proposed algorithm to solve the FVM-discretized Richards equation using a message passing mechanism.}
    \label{fig_flowchart}
\end{figure}

Once the adaptive fixed-point iteration scheme for the FVM-discretized Richards equation is established, we incorporate it in our MP-FVM algorithm to enhance the solver accuracy and ability to retain underlying physics (e.g., mass conservation). As discussed previously, the message passing neural PDE solver proposed by \citet{brandstetter} comprises three main components: an encoder, a processor, and a decoder. The message passing mechanism is implemented within the processor that operates in the latent space. However, it has not been extended to discretized PDEs. In this work, we introduce the message passing mechanism for the discretized Richards equation by defining a latent variable $\mu_i^{m,s}$ as the processor. Therefore, by leveraging our adaptive fixed-point iteration scheme, we can now solve the latent variable iteratively to enhance the convergence and numerical stability of the message passing mechanism. Specifically, this integrative algorithm, MP-FVM, adopts one neural network (encoder) $\hat{f}_{\text{NN}}$ to learn the map $\psi_i^{m,s} \mapsto \mu_i^{m,s}$ and another neural network (decoder) $\hat{f}_{\text{NN}}^{-1}$ to learn the inverse map $\mu_i^{m,s} \mapsto \psi_i^{m,s}$. Overall, our MP-FVM algorithm involves offline training (dataset preparation and encoder-decoder training) and solution (message passing) process, which are summarized in the flowchart of Figure \ref{fig_flowchart}.

\subsection{Dataset Preparation and Data Augmentation}
The dataset used to train the encoder and decoder neural networks comes from two different sources/solvers. Specifically, for each cell $V_i$ and time step $m$, we approximate the latent variable solution $\mu_{i}^{m,S}$ from a finite difference solver (e.g., \citet{ireson2023simple}). Here, $S$ is the user-specified total iteration number. The corresponding $\psi_{i}^{m,S}$ solution is obtained separately from the fixed-point iteration scheme of Equation \eqref{eqn_lscheme2} using a static parameter $\tau$. The resulting set of solution pairs, $\big\{\big(\psi_{i}^{m,S},\mu_{i}^{m,S}\big)\big\}_{i,m}$, form a set of original ``reference solutions''. In actual implementation, we obtain multiple sets of original reference solutions by selecting multiple total iteration numbers ($S_1,\dots,S_p$) and/or fixed-point parameters ($\tau_1,\dots,\tau_r$) that cover their ranges expected during the actual solution process. These sets of original reference solutions, which are $\big\{\big(\psi_{i}^{m,S_1},\mu_{i}^{m,S_1}\big)|_{\tau_1}\big\}_{i,m}, \dots, \big\{\big(\psi_{i}^{m,S_p},\mu_{i}^{m,S_p}\big)|_{\tau_r}\big\}_{i,m}$, are combined to form a larger set to perform data augmentation. 

Next, to apply data augmentation, we introduce Gaussian noise $Z_q\sim \mathcal{N}(0,\,\sigma_q^{2})$ with different variances $\sigma_1^2,\dots, \sigma_Q^2$ to each and every element in the reference solution set obtained previously. After data augmentation, the resulting expanded set of reference solutions, $\big\{(\psi_{i}^{m,S_1} + Z_p,\mu_{i}^{m,S_1} + Z_q)|_{\tau_1}\big\}_{i,m,q}, \dots, \big\{(\psi_{i}^{m,S_p} + Z_q,\mu_{i}^{m,S_p} + Z_q)|_{\tau_r}\big\}_{i,m,q}$, is denoted as $\mathcal{S}$ and will be used for neural network training. This data augmentation step not only increases the size of the training dataset, but also reflects the characteristics of actual soil sensing data, which are subject to various measurement uncertainties. Furthermore, In Section \ref{sec_1d}, we will show that introducing Gaussian noise can greatly reduce the biases of reference solutions and enhance generalization performance \citep{da2011pca}, thereby significantly improving the accuracy of numerical solutions.

\subsection{Neural Network Training}
A neural network is capable of approximating any function provided that it contains enough neurons \citep{hornik1991approximation,pinkus}. In the actual implementation, depending on the problem settings, the desired choices of optimal optimizer, number of hidden layers, and activation functions can vary. Based on our extensive research and hyperparameter tuning, we find that a simple three-layer neural network with 256 neurons in each layer achieves the best performance for most 1-D through 3-D problems compared to other more complex neural network architectures (e.g., LSTM). Also, we find that stochastic gradient decent (SGD) optimizer often outperforms others (e.g., Adam or RMSProp). The learning rate is set to be $0.001$. This simple neural network structure makes our MP-FVM algorithm training much less computationally expensive compared to state-of-the-art neural PDE solvers (e.g., \citet{lu2020,brandstetter}). 

In terms of loss function design, we note that the solution of the Richards equation at a given time step depends on the pressure head solution at the initial condition and previous time steps. A small perturbation in these solutions can lead to slow convergence or inaccurate solutions at the final time step. To account for this, we introduce Sobolev training \citep{czarnecki2017} for both neural networks $\hat{f}_{\text{NN}}$ and $\hat{f}^{-1}_{\text{NN}}$ to ensure compatibility and stability in the same solution space. We implement Sobolev training by adding a Sobolev regularization term to the standard Mean Squared Error (MSE) in the loss functions for $\hat{f}_{\text{NN}}$ and $\hat{f}^{-1}_{\text{NN}}$:
\begin{equation}\label{eqn_sobolev}
   \mathcal{L}_{\hat{f}_{\text{NN}}} = \underbrace{\frac{1}{|\mathcal{S}|} \sum_{(\psi, \mu) \in \mathcal{S}} \left( \mu - \hat{f}_{\text{NN}}(\psi) \right)^2}_{\text{MSE term}} + \lambda_{\hat{f}_{\text{NN}}} \cdot \underbrace{\frac{1}{|\mathcal{S}|} \sum_{(\psi, \mu) \in \mathcal{S}} \left(\left\| \nabla \left(\mu - \hat{f}_{\text{NN}}(\psi)\right) \right\|_{L^2}^2 \right)}_{\text{Sobolev regularization term}},
\end{equation}
and 
\begin{equation}\label{eqn_sobolev_2}
   \mathcal{L}_{\hat{f}^{-1}_{\text{NN}}} = \underbrace{\frac{1}{|\mathcal{S}|} \sum_{(\psi, \mu) \in \mathcal{S}} \left( \psi - \hat{f}^{-1}_{\text{NN}}(\mu) \right)^2}_{\text{MSE term}} + \lambda_{\hat{f}^{-1}_{\text{NN}}} \cdot \underbrace{\frac{1}{|\mathcal{S}|} \sum_{(\psi, \mu) \in \mathcal{S}} \left( \left\| \nabla \left(\psi - \hat{f}^{-1}_{\text{NN}}(\mu)\right) \right\|_{L^2}^2 \right)}_{\text{Sobolev regularization term}},
\end{equation}
where $\lambda_{\hat{f}_{\text{NN}}}$ and $ \lambda_{\hat{f}^{-1}_{\text{NN}}}$ are user-specified regularization parameters for the neural networks $ \lambda_{\hat{f}^{-1}_{\text{NN}}}$ and $\hat{f}^{-1}_{\text{NN}}$, respectively. Here, we use the Leaky ReLU activation function, as it has been shown that there exists a single hidden-layer neural network with ReLU (or Leaky ReLU) activation function that can approximate any function in a Sobolev space \citep{czarnecki2017}. Overall, this combined loss function ensures that the model not only produces accurate predictions but also generates smooth and regular outputs by matching the gradients of the true function.

\subsection{Message Passing Process}
When neural network training is complete, the trained encoder $\hat{f}_{\text{NN}}$ and decoder  $\hat{f}^{-1}_{\text{NN}}$ can then be incorporated into Equation \eqref{eqn_lscheme2} to derive the following fixed-point iterative scheme for the latent variables with message passing mechanism: 
\begin{equation}\label{eqn_DRWlscheme}
     \mu_i^{m+1,s+1} = \mu_i^{m+1,s}+\tau_i^{m+1,s}\sum_{j\in \mathcal{N}_i}K_{\omega_{i,j}}^{m+1,s}\mathbf{e} \cdot \mathbf{n}_{\omega_{i,j}}\frac{\mu_j^{m+1,s}-\mu_i^{m+1,s}}{\text{dist}(V_j,V_i)}A_{\omega_{i,j}}+ \hat{f}_{\text{NN}}(J),
\end{equation}
where $J = \tau_i^{m+1,s} \sum_{j\in \mathcal{N}_i}K_{\omega_{i,j}}^{m+1,s} \mathbf{e} \cdot \mathbf{n}_{\omega_{i,j}}\frac{z_j-z_i}{\text{dist}(V_j,V_i)}A_{\omega_{i,j}}-\tau_i^{m+1,s}\left(\frac{\theta_i^{m+1,s}-\theta_i^m}{\Delta t} + S(\psi_i^{m+1,s})\right)\text{vol}(V_i)$. To solve Equation \eqref{eqn_DRWlscheme}, we will adopt a similar strategy as in Equation \eqref{eqn_Lbound2} to adaptively select the linearization parameter $\tau_i^{m+1,s}$. To start the message passing process, we obtain the initial pressure head solutions in the control volume at $m=0$ from the initial and boundary conditions. These initial pressure head solutions can be mapped to the latent space via trained encoder network $\hat{f}_{\mathrm{NN}}$. Next, for each new time step $m+1$, the latent variable for every cell can be iteratively solved by Equation \eqref{eqn_DRWlscheme} by utilizing the trained neural networks  $\hat{f}_{\mathrm{NN}}$ and $\hat{f}^{-1}_{\mathrm{NN}}$. Note that the iterative usage of $\hat{f}^{-1}_{\mathrm{NN}}$ is implicitly implied in the MP-FVM algorithm, as the term $J$ in Equation \eqref{eqn_DRWlscheme} contains $\psi_i^{m+1,s}$ that must be evaluated by applying $\hat{f}^{-1}_{\mathrm{NN}}$ on latent variable $\mu_i^{m+1,s}$. Also, it is worth mentioning that, since $\psi$ and $J$ have different scales, in actual implementation, in addition to $\hat{f}_{\mathrm{NN}}$ for learning $\psi_i^{m,s} \rightarrow \mu_i^{m,s}$, we train another neural network named $\hat{f'}_{\text{NN}}$ for mapping $J$ to the latent space in Equation \eqref{eqn_DRWlscheme}. To monitor convergence of the iterative message passing process, we define the relative error $\mathrm{RE}_s$ as:
\begin{equation}\label{eqn_error}
    \mathrm{RE}_s \coloneqq \frac{||\mu^{m+1,s+1} - \mu^{m+1,s}||_{L^2}}{||\mu^{m+1,s+1}||_{L^2}},
\end{equation}
where $\mu^{m+1,s+1} = (\mu_1^{m+1,s+1},\dots, \mu_N^{m+1,s+1})^T$ and so on. Once $\mathrm{RE}_s$ is below a user-specified tolerance $\mathrm{tol}$ (typically in the order of $10^{-6}$), we declare convergence of $\{\mu_i^{m+1,s}\}_s$ to $\mu_i^{m+1}$. From there, one can determine the converged $\psi_i^{m+1}$ using $\hat{f}^{-1}_{\text{NN}}$, followed by obtaining other physical quantities such as soil moisture content $\theta_i^{m+1}$ and $\textbf{q}_i^{m+1}$ from the WRC and HCF models (Table \ref{table_hcfwrc}) and Equation \eqref{eqn_richards}. The entire solution process then repeats itself in the next time step until $m = \lceil \frac{T}{\Delta t}\rceil -1$.

Furthermore, it is worth mentioning that, when neural network training for a specific problem setting (e.g., boundary condition and initial condition) is complete, the trained neural networks can be saved as a pre-trained model. As we encounter a new problem setting, the pre-trained model provides a strong starting point that can be quickly refined with a small number of epochs (typically no more than 100) before it can be deployed to solve the new problem. The use of pre-trained model is a well-established technique in machine/deep learning for leveraging knowledge learned from (large) datasets, reducing the need for extensive training data and computation, and enabling faster deployment and improved performance in new tasks through fine-tuning.

\subsection{Convergence of MP-FVM Algorithmn}

The convergence of our MP-FVM algorithm, which features the sequence $\{\mu^{m+1,s}\}_s$, can be established by extending Theorem \ref{theorem_convergence} and investigating the convergence behavior of stochastic gradient descent (SGD) for neural network realizations of $\hat{f}_{\mathrm{NN}}$ and $\hat{f}^{-1}_{\mathrm{NN}}$. Similar to Theorem \ref{theorem_convergence}, we consider functions $\{\mu^{m+1,s}\}_s$ and $\mu^{m+1}$ instead of their discretized variants.
\begin{theorem}\label{theorem_convergence2}
    The sequence $\{\mu^{m+1,s}\}_s$ converges to $\mu^{m+1}$ for $m = 0,1,\ldots,\lceil\frac{T}{\Delta t}\rceil-1$.
\end{theorem}
\begin{proof}
    See Appendix \ref{appen_c} for the complete proof.
\end{proof}

\section{Case Studies} \label{sec_casestudies}

Now that we have introduced the MP-FVM algorithm formulation for the Richards equation, in this section, we evaluate our MP-FVM framework on a series of 1-D through 3-D benchmark problems modified from the literature \citep{celia1990general,gkasiorowski2020numerical,tracy2006clean,berardi20181d,orouskhani_impact_2023}. Specifically, we extensively study the 1-D benchmark problem of \citet{celia1990general} to demonstrate the need and benefits of different components employed in our MP-FVM algorithm, including adaptive fixed-point iteration scheme, encoder-decoder architecture and message passing mechanism, and Sobolev training. Also, using this problem as a benchmark, we demonstrate the accuracy of our solution algorithm with respect to state-of-the-art solvers. In the 1-D layered soil case study proposed by \citet{berardi20181d}, we show that our MP-FVM algorithm is capable of handling discontinuities in soil properties and modeling the infiltration process through the interface of two different soils. In the 2-D case study adopted from \citet{gkasiorowski2020numerical}, we show that our MP-FVM algorithm can better satisfy the mass balance embedded in the Richards equation. In the 3-D case study adopted from \citet{tracy2006clean} in which an analytical solution to the Richards equation exists, we show that our MP-FVM algorithm produces much more accurate solutions compared to conventional FVM solvers. Finally, we study a 3-D problem adopted from \citet{orouskhani_impact_2023} featuring an actual center-pivot system and validate the accuracy and robustness of our MP-FVM algorithm in modeling real-world precipitation and irrigation scenarios for a long period of time.

\subsection{A 1-D Benchmark Problem} \label{sec_1d}
Here, we study the 1-D benchmark problem over a $40$ cm deep soil presented by \citet{celia1990general}. The HCF and WRC adopt the model of \citet{haverkamp1977comparison} (see Table \ref{table_hcfwrc}), whose parameters are listed in Table \ref{table_1dparameters}. The initial condition is given by $\psi(z,0)=-61.5$ cm, whereas the two boundary conditions are $\psi(40\text{ cm},t)=-20.7\text{ cm}, \, \psi(0,t)=-61.5$ cm, respectively \citep{haverkamp1977comparison}. This benchmark problem ignores the sink term.

\begin{table}[ht!]\label{table_1dparameters}
    \centering
    \begin{tabular}{ccc}
    \toprule
    Soil-specific Parameters & Values & Units\\
    \midrule
    Saturated hydraulic conductivity, $K_s$ & $0.00944$ & cm/s\\
    Saturated soil moisture content, $\theta_s$ & $0.287$ & --\\
    Residual soil moisture content, $\theta_{r}$ & $0.075$ & --\\
    $\alpha$ in Haverkamp's model & $1.611\times 10^6$& cm\\
    $A$ in Haverkamp's model & $1.175\times 10^6$ & cm \\
    $\beta$ in Haverkamp's model & $3.96$& --\\
    $\gamma$ in Haverkamp's model & $4.74$ & --\\
    Total time, $T$ & $360$ & s\\
    \bottomrule
    \end{tabular}
    \caption{soil-specific parameters and their values used in the 1-D case study of \citet{celia1990general} based on the empirical model developed by \citet{haverkamp1977comparison}.}
\end{table}

Through this 1-D illustrative example, we will highlight the benefits of (a) adopting an adaptive fixed-point iteration scheme as opposed to standard the fixed-point iteration scheme, (b) implementing the MP-FVM algorithm as opposed to the conventional FVM method, and (c) integrating the adaptive fixed-point iteration scheme with encoder-decoder network and message passing mechanism in a holistic numerical framework.

\subsubsection{The Need for Adaptive fixed-point iteration scheme}

\begin{figure}[ht!]
    \centering
    \includegraphics[trim={2cm 0 3cm 0},clip,width=\textwidth]{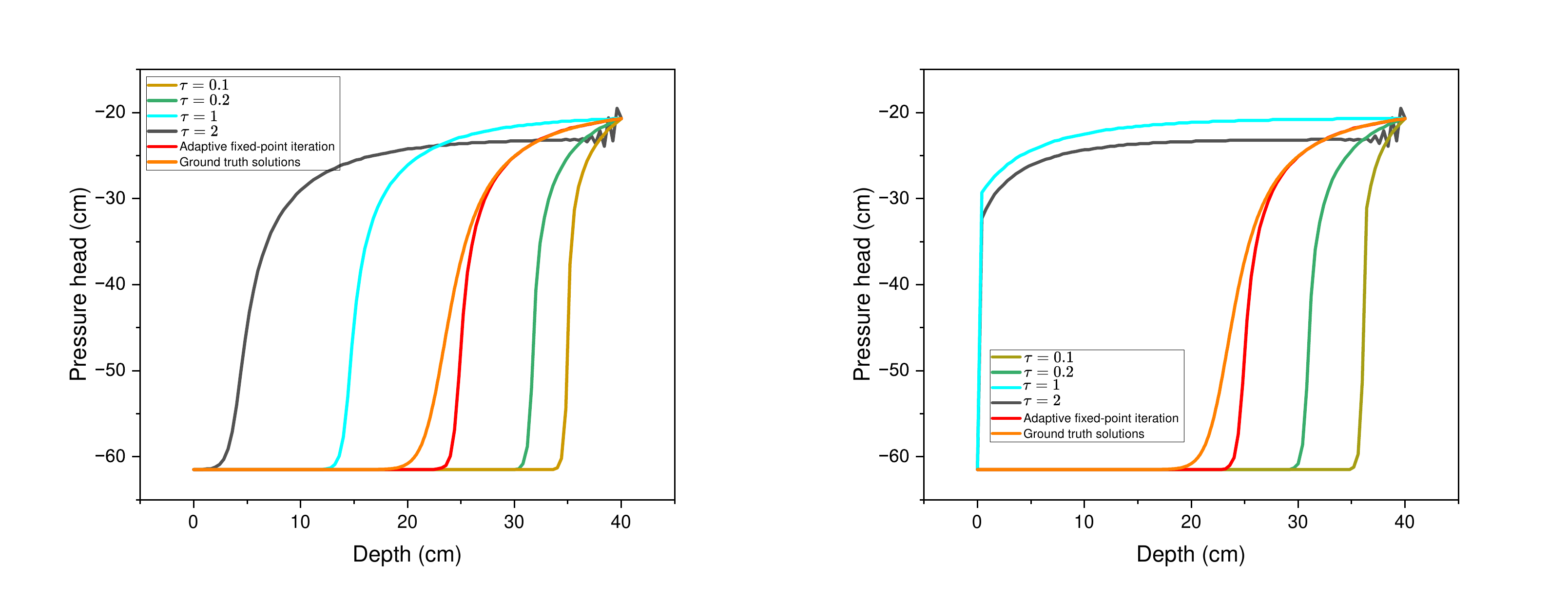}
    \vspace{-3em}
    \caption{Comparison of pressure head solution profiles at $t=T=360$ seconds under (a) $S = 500$ iterations and (b) $\mathrm{tol}=3.2\times 10^{-5}$ for the 1-D benchmark problem \citep{celia1990general} using standard and adaptive fixed-point iteration schemes (Equation \eqref{eqn_lscheme2}). The solutions obtained from \citet{celia1990general} based on very fine space and time steps are marked as the ground truth solutions.}
    \label{fig_1dLscheme}
\end{figure}

To illustrate how adaptive fixed-point iteration scheme improves convergence and accuracy of conventional fixed-point iteration schemes, we compare the pressure head solution profiles at $t = T = 360$ seconds obtained by different static fixed-point parameters after (a) $S=500$ iterations and (b) $\mathrm{tol}=3.2\times 10^{-5}$. We adopt a spatial grid containing $101$ mesh points ($\Delta z = 0.4$ cm) and a temporal grid satisfying the Courant-Friedrichs-Lewy (CFL)-like condition, typically expressed as $\Delta t \leq \frac{\Delta z^2}{2K}$ \citep{de2013courant}. As shown in Figure \ref{fig_1dLscheme}, when using static fixed-point iteration scheme, the choice of parameter $\tau$ and the total number of iterations can impact the solution accuracy and algorithm stability significantly. For example, when the fixed-point parameter is too large (e.g., $\tau=2$ for this problem), the stability of the static fixed-point iteration scheme can be adversely affected (as illustrated by the zigzag pressure head profile towards $z = 40$ cm). Another key observation is that, increasing the total number of iterations sometimes deteriorates solution accuracy of static fixed-point iteration scheme. These observations pose practical challenges for using static fixed-point iteration scheme, especially when the ground truth solutions are absent, as identifying the optimal fixed-point parameter and total number of iterations that would yield accurate solutions will not be possible without referring to ground truth solutions. This motivates us to develop adaptive fixed-point iteration scheme as a robust and reliable numerical scheme that produces solutions that are close to ground truth solutions without trail-and-error parameter tuning. Also, it is worth noting that our adaptive fixed-point iteration scheme successfully bypasses the singularity issue as $\frac{1}{\tau_i^{m+1,s}}$ approaches to 0 and correctly calculates the pressure head solutions for $z \in [0,20 \text{ cm}]$ where $\dot \theta(\psi)$ becomes small.

\subsubsection{The Need for Encoder-Decoder Architecture}

To generate the reference solutions, we consider a coarse spatial discretization containing $40$ cells (i.e., grid size $\Delta z = 1$ cm) and solve for $T = 360$ seconds. The time step size $\Delta t$ is determined using the CFL condition \citep{de2013courant}. A set of pressure head solutions $\psi$ is obtained using the finite difference method that incorporates a modified Picard iteration scheme developed by \citet{celia1990general}. Meanwhile, another set of pressure head solutions, which essentially becomes the latent variable dataset $\mu$ for neural network training, is obtained from the fixed-point iteration scheme of Equation \eqref{eqn_lscheme2} under $4$ different static fixed-point parameter $\tau=0.25,\, 0.24,\, 0.23,\, 0.22$ and $10$ different total iteration counts $S = 1,000$, $2,000$, up to $10,000$. 

\begin{figure}[ht!]
    \centering
    \includegraphics[trim={0 0 0 4em}, clip, width=0.9\textwidth]{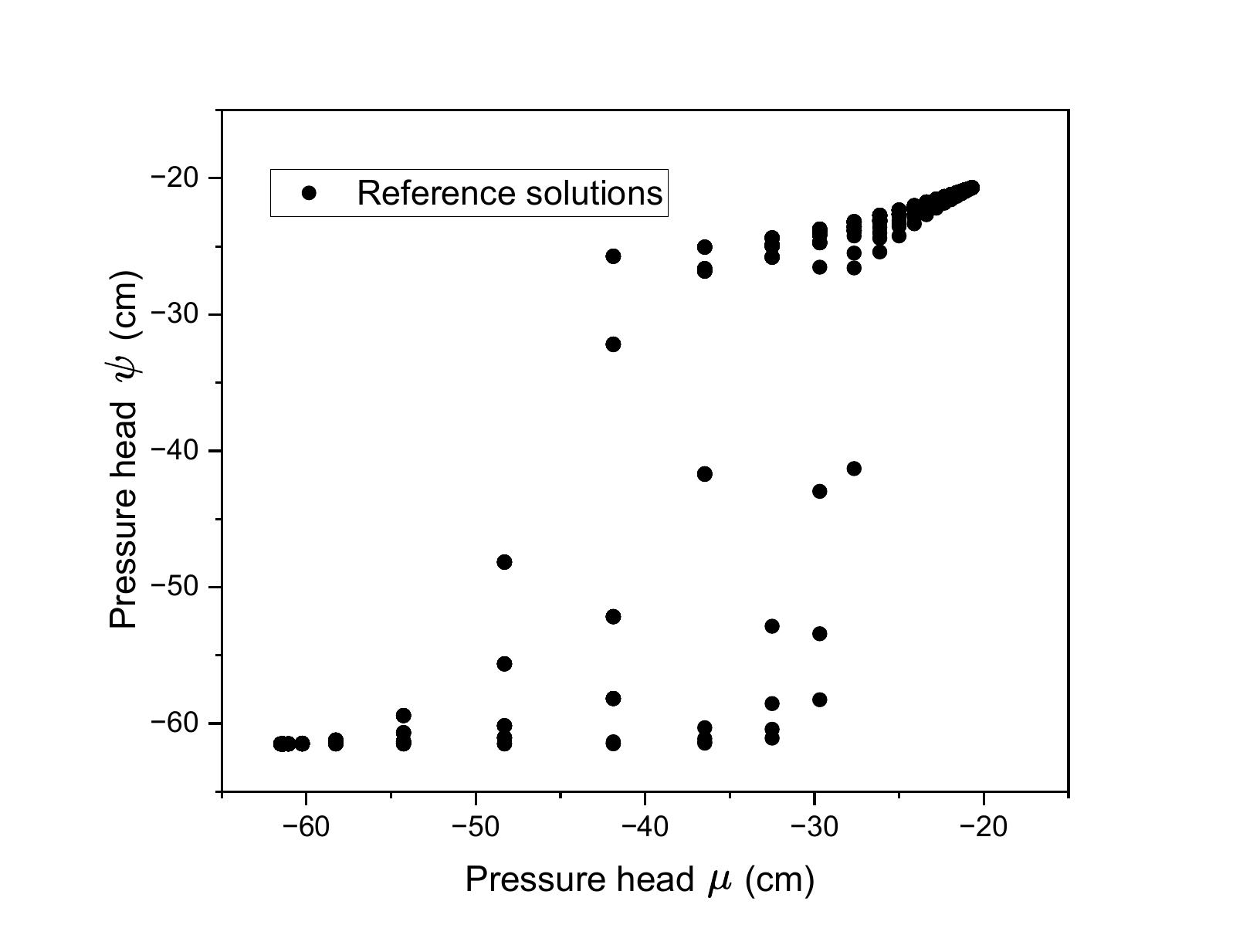}
    \vspace{-2em}
    \caption{The relationships between $1640$ pressure head solutions $\psi$ and $\mu$, which are obtained by two distinct approaches. The resulting nonlinearity present in these reference solutions highlights need for data-driven approach.}
    \label{fig_assumption}
\end{figure}

As mentioned earlier, reference solutions utilized to train the encoder $\hat{f}_{\text{NN}}$ and decoder  $\hat{f}^{-1}_{\text{NN}}$ come from two different sources. As shown in Figure \ref{fig_assumption}, a highly nonlinear relationship between two sources of pressure head solutions is observed. This is mainly because pressure head solutions from different sources exhibit different sensitivities with respect to different choices of $\tau$ and $S$. Without knowing the ground truth solutions a priori, it is hard to determine which set of pressure head solutions is more accurate. This motivates us to adopt an encoder-decoder architecture to explicitly capture this nonlinear relationship, which encapsulates the sensitivity of solution with respect to different choices of $\tau$ and $S$. 

\begin{figure}[ht!]
    \centering
    \includegraphics[trim={0 0 0 0em}, clip, width=\textwidth]{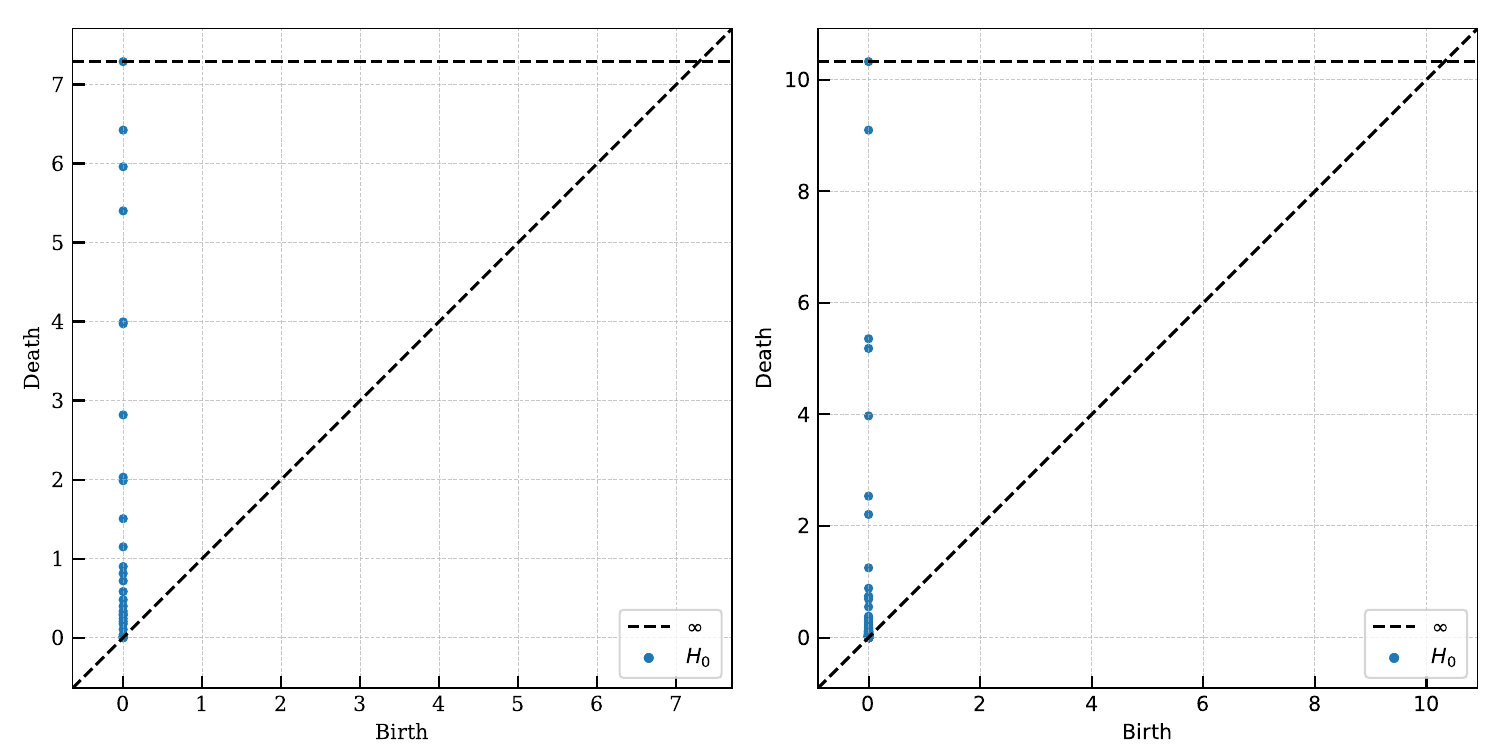}
    \vspace{-2em}
    \caption{Persistence diagrams \citep{edelsbrunner2013persistent} for pressure head solutions $\psi$ (left) and $\mu$ (right). The marked differences in topological features illustrate the need for an encoder to map $\psi$ into the topological space of $\mu$. Here, $\infty$ refers to infinite lifespan and $H_0$ are connected components.}
    \label{fig_ph}
\end{figure}

Another motivation for adopting an encoder-decoder architecture in our numerical solver comes from the fact that different sources of pressure head solutions also exhibit different topological features. To see this, we use persistent homology \citep{edelsbrunner2013persistent} as a way to capture the multiscale topology of each source of pressure head solutions. Specifically, we construct a sequence of simplicial complexes and track the ``birth'' and ``death'' of topological features across this sequence. Figure \ref{fig_ph} shows that the $\mu$ solutions exhibit longer-lasting topological components than the $\psi$ solutions, as all points die off much sooner (e.g., $\sim 7.4$ on the death axis) for the $\psi$ solutions. Therefore, the use of an encoder $\hat{f}_{\text{NN}}$, which maps the pressure head solutions $\psi$ to a latent space where $\mu$ solutions lie, can capture the distinct topological structures of two sources of pressure head solutions. Similarly, the decoder $\hat{f}^{-1}_{\text{NN}}$ transforms the latent representation $\mu$ back to the original solution space, ensuring that the essential topological features of $\psi$ solutions are accurately captured and reconstructed.

\subsubsection{Improving MP-FVM Algorithm Performance via Sobolev Training and Encoder-decoder Architecture}

\begin{figure}[ht!]
    \centering
    \includegraphics[trim={0 0 0 4em}, clip, width=\textwidth]{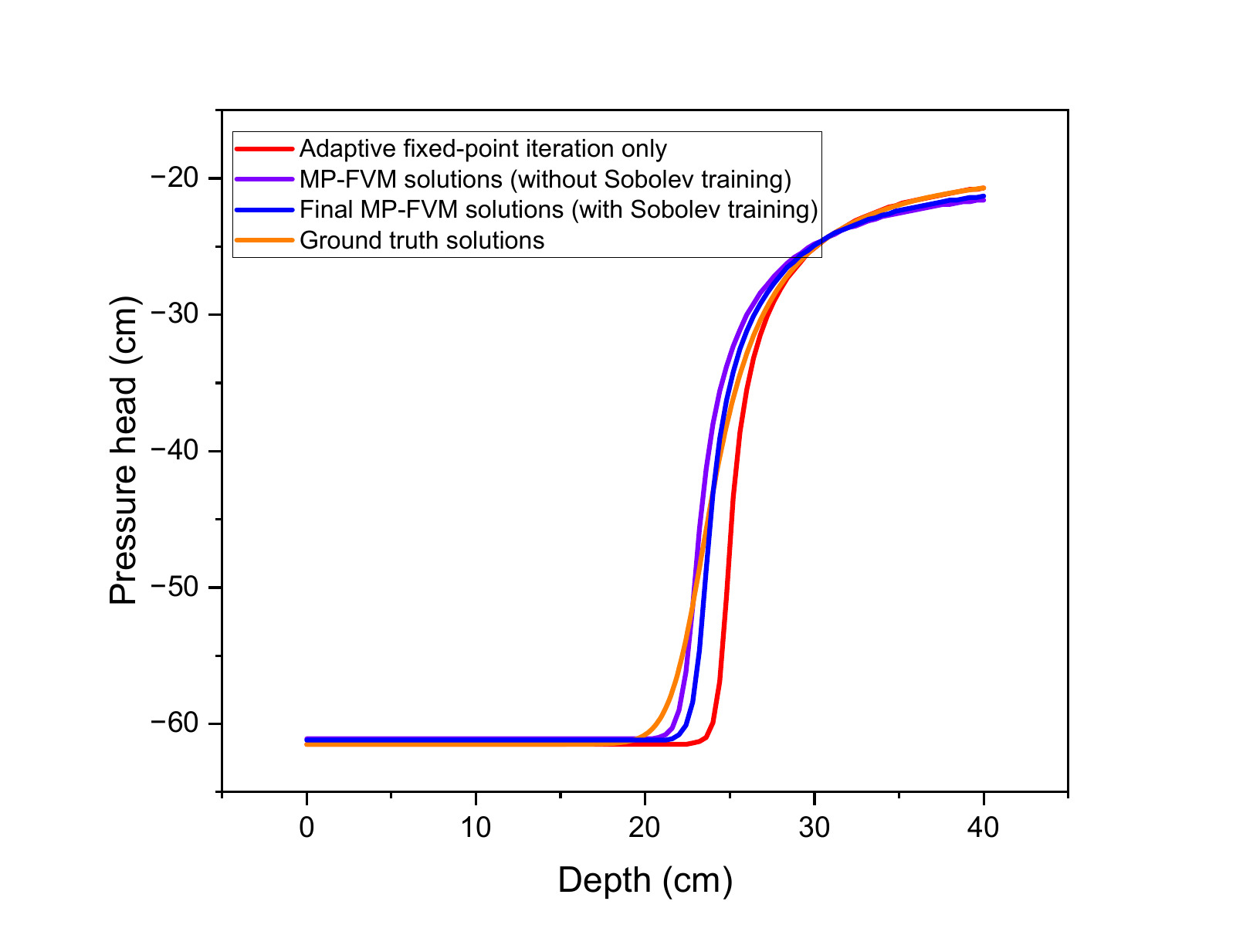}
    \vspace{-1em}
    \caption{Comparison of pressure head solution profiles at $t=T=360$ seconds produced from adaptive fixed-point iteration scheme only (Equation \eqref{eqn_lscheme2}) and from MP-FVM algorithm (Equation \eqref{eqn_DRWlscheme}) with and without implementing Sobolev training.}\label{fig_integration}
\end{figure}

As previously discussed, we perform data augmentation on the reference solutions to increase dataset size and enhance generalization performance. Specifically, after we obtain a set of $\psi$ solutions using the finite difference method developed by \citet{celia1990general}, we make multiple copies of it and append each copy to the $\mu$ solutions obtained by the fixed-point iteration scheme of Equation \eqref{eqn_lscheme2} under different static $\tau$ and $S$ values. We then add zero-mean Gaussian noises with standard deviation varying from $0.1$ to $0.5$ to these augmented reference solutions. Overall, this leads to a total of $17,097$ reference solutions for neural network training and validation. Note that, as previously discussed, the original and augmented reference solutions are generated using a coarse grid ($\Delta z = 1$ cm). Thus, they can be obtained relatively efficiently. On the other hand, in the solution step, we will use a more refined grid containing $101$ mesh points ($\Delta z = 0.4$ cm). This ``coarse-to-fine'' approach can therefore enhance the solution accuracy of our MP-FVM algorithm without requiring a large amount of high-accuracy, fine-mesh training data. Furthermore, when augmented reference solutions are used for training, only $100$ additional epochs are needed to retrain neural networks that have already been trained using the original reference solutions. Second, we notice that there is only a slight difference in the final pressure head solution profile when Gaussian noises of different magnitudes are directly added to the original reference solutions without augmenting them together. Third, increasing training data size (from $1,640$ to $17,097$) via data augmentation of original reference solutions is an effective way to improve solution accuracy of our MP-FVM algorithm, as the pressure head profile matches very well with the ground truth solution.

From Figure \ref{fig_integration}, it is clear that integrating adaptive fixed-point iteration scheme in the MP-FVM framework synergistically improves the overall solution accuracy of the Richards equation, especially in the region where pressure head changes rapidly with respect to depth (i.e., between $z=20$ to 30 cm). On the other hand, we observe slight discrepancy in pressure head solution close to $z=40$ cm when comparing our MP-FVM algorithm with ground-truth solutions, whereas the solution produced by adaptive fixed-point iteration scheme alone matches perfectly with ground-truth solution at $z=40$ cm, which corresponds to one of the boundary conditions. We believe that this is due to the fact that $\hat{f}_{\text{NN}}$ and $\hat{f}^{-1}_{\text{NN}}$ only approximate the true relationships $f$ and $f^{-1}$, respectively, and the resulting induced error causes discrepancies in pressure head solutions even at the boundaries. To overcome this limitation, one way is to increase the size of the augmented reference solutions for neural network training. Another approach is to switch from MP-FVM (Equation \eqref{eqn_DRWlscheme}) to adaptive fixed-point iteration scheme only (i.e., Equation \eqref{eqn_lscheme2}) when solving for the boundary conditions. We leave this refinement for future research.

\begin{figure}[ht!]
    \centering
    \includegraphics[trim={2cm 0.5cm 3cm 2cm},clip,width=0.85\textwidth]{fig_lambda_s2.jpg}
    \vspace{-1em}
    \caption{Comparison of pressure head solution profiles at $t=T=360$ seconds produced from MP-FVM algorithm (Equation \eqref{eqn_DRWlscheme}) with implementing Sobolev training with different regularization parameters in Equation \eqref{eqn_sobolev} and Equation \eqref{eqn_sobolev_2} at Scenario 2. Here, we use the same $\lambda=\lambda_{\hat{f}_{\text{NN}}}=\lambda_{\hat{f}^{-1}_{\text{NN}}}$ and all neural networks are trained from scratch.}\label{fig_sobolev}
\end{figure}

Figure \ref{fig_sobolev} illustrates how Sobolev training affects the solution quality of our MP-FVM algorithm. Specifically, we find that, first, the effectiveness of Sobolev training depends on the choice of hyperparameter $\lambda$. Second, larger values of $\lambda$ (e.g., $10^{-5}$) may not lead to improved accuracy in pressure head solution, as in this case, neural network training may prioritize smoothness or derivative agreement over fitting the pressure head solutions. Third, smaller values of $\lambda$ (e.g., $10^{-9}$) could still be useful in improving solution accuracy compared to without Sobolev training (i.e., $\lambda=0$). Last but not least, we notice that, when pre-trained models are used, the sensitivity of pressure head solution to $\lambda$ is significantly reduced, especially for $\lambda < 10^{-6}$. We suspect that this is because pre-trained models already capture the relationships between $\psi$ and $\mu$ solutions reasonably well, so that the Sobolev loss primarily serves to fine tune the models.

\subsubsection{Convergence and Solution Accuracy Comparison} 

We compare our MP-FVM algorithm with other solvers based on computational performance and solution accuracy under two scenarios. In Scenario 1, we set the error tolerance $\mathrm{tol}$ to be $3.2\times10^{-5}$, whereas in Scenario 2, we set the total number of iterations $S = 500$. For static fixed-point iteration scheme, we use an optimal fixed-point parameter $\tau = \frac{1}{3.5} \approx 0.2857$ identified by trail-and-error process. In terms of computational performance, we use the condition number of matrix $\mathbf{A}$ defined in Equation \eqref{eqn_matrix}, which measures the sensitivity of fixed-point iteration scheme subject to small perturbations, as the metric. 

\begin{table}[ht!]
    \centering
    \begin{tabular}{ccc}
        \toprule
        \multirow{2}{*}{Algorithm} & \multicolumn{2}{c}{Average condition number of $\mathbf{A}$ obtained from \citet{zarba1988thesis} (Scenario 1)} \\ \cline{2-3} 
        & Static fixed-point iteration scheme & Adaptive fixed-point iteration scheme \\ \midrule
        FVM & $1.7668$    & $1.0064$ \\ 
        MP-FVM & $1.7419$     & $1.0075$\\ \bottomrule
    \end{tabular}
    \caption{Comparison of average condition number under Scenario 1 across all time steps (as Equation \eqref{eqn_matrix} already considers all discretized cells) for conventional FVM and our MP-FVM algorithms that implement static or adaptive fixed-point iteration scheme.}\label{table_1dcon1}
\end{table}

\begin{table}[ht!]
    \centering
    \begin{tabular}{ccc}
        \toprule
        \multirow{2}{*}{Algorithm} & \multicolumn{2}{c}{Average condition number of $\mathbf{A}$ obtained from \citet{zarba1988thesis} (Scenario 2)} \\ \cline{2-3} 
        & Static fixed-point iteration scheme & Adaptive fixed-point iteration scheme \\ \midrule
        FVM & $1.7206$    & $1.0064$ \\ 
        MP-FVM & $1.7113$     & $1.0071$\\ \bottomrule
    \end{tabular}
    \caption{Comparison of average condition number under Scenario 2 across all time steps for conventional FVM and our MP-FVM algorithms that implement static or adaptive fixed-point iteration scheme.}\label{table_1dcon2}
\end{table}

From Tables \ref{table_1dcon1} and Table \ref{table_1dcon2}, we see that implementing adaptive fixed-point iteration scheme significantly improves the stability of conventional FVM and our MP-FVM algorithms, as matrix $\mathbf{A}$ is well-conditioned. These observations suggest that adaptive fixed-point iteration scheme outperforms static fixed-point iteration scheme in enhancing the convergence behavior of discretization-based solvers.

In terms of solution accuracy, we consider two metrics. The first metric is the discrepancy from the ground truth solutions of \citet{celia1990general}. The comparison results are illustrated in Figure \ref{fig_1dcompare}. The second metric is the solver's performance in preserving the mass (moisture) balance, which is quantified by the mass balance measure $\mathrm{MB}$ defined in \citet{celia1990general}:
\begin{equation} \label{eqn_mb}
    \mathrm{MB} =\frac{\text{total additional mass in the domain}}{\text{total water flux into the domain}}.
\end{equation}

In Figure \ref{fig_1dcompare}, we compare the pressure head profiles obtained from our MP-FVM algorithm (which implements adaptive fixed-point iteration scheme and Sobolev training), the conventional FVM algorithm (that implements adaptive fixed-point iteration scheme), and a state-of-the-art physics-informed neural network (PINN) solver based on \citet{pinn}, against the ground truth solution \citep{celia1990general}. Clearly, in both scenarios, compared with the MP-FVM solutions, PINN and FVM solutions are further apart from ground truth solutions.

\begin{figure}[ht!]
    \centering
    \includegraphics[width=\textwidth]{fig_1dcompare.jpg}
    \vspace{-2em}
    \caption{Pressure head profiles at $t=T=360$ sec obtained by different algorithms under (left) Scenario 1, and (right) Scenario 2. Both conventional FVM and our MP-FVM algorithm incorporate adaptive fixed-point iteration scheme. Note the PINN solver is not an iterative method, thus the solution profile is the same under both scenarios.} \label{fig_1dcompare}
\end{figure}

From Tables \ref{table_1dMB1} and \ref{table_1dMB2}, we observe that, in both Scenarios 1 and 2, our MP-FVM algorithm achieves the best $\mathrm{MB}$ values when using either coarse time steps suggested by the CFL condition \citep{de2013courant} or a fixed time step. Considering that using coarse time steps reduces solution time without affecting solution quality, adopting a CFL-like condition is desired.

\begin{table}[ht!]
    \centering
    \begin{tabular}{llll}
    \toprule
        Method used & Scenario & Average $\Delta t\text{ (sec)}$& $\mathrm{MB}$  \\ \midrule
        FVM algorithm& 1 & $18.90$ &$96.13\%$\\
        MP-FVM algorithm & 1 & $18.68$ & \pmb{$100.23\%$}  \\ 
        FVM algorithm & 2 & $17.62$ & $86.04\%$  \\
        MP-FVM algorithm & 2 & $18.35$ &\pmb{$97.29\%$}  \\ 
        \citet{celia1990general} & N/A & $10$ & $95.00\%$\\
        \bottomrule
    \end{tabular}
    \caption{$\mathrm{MB}$ results of different numerical methods. Note that here, $\Delta t$ is the determined for each method by the CFL condition \citep{de2013courant} and we take the average across all iterations.} \label{table_1dMB1}
\end{table}

\begin{table}[ht!]
    \centering
    \begin{tabular}{llll}
    \toprule
        Method used & Scenario & $\mathrm{MB}$ ($\Delta t =15$ sec)  \\ \midrule
        FVM algorithm & 1 & $98.87\%$ \\
        MP-FVM algorithm & 1 &  \pmb{$100.72\%$}  \\ 
        FVM algorithm & 2 & $96.79\%$ \\
        MP-FVM algorithm & 2 & \pmb{$97.81\%$}  \\ 
        \citet{celia1990general} & N/A & $95.00\%$\\
        \bottomrule
    \end{tabular}
    \caption{$\mathrm{MB}$ results of different numerical methods, in which a common $\Delta t = 10$ seconds is used for all numerical methods.} \label{table_1dMB2}
\end{table}

\subsubsection{Remark on Computational Efficiency}
Although our MP-FVM framework does involve neural network training which will take some additional time, there are several well-established strategies widely used in the machine/deep learning community to reduce the overall computational time and costs. For example, as previously discussed, one can leverage the previously trained neural network from a different problem setting as a good starting point to train with new dataset for the new problem setting in just a small number of epochs. To see this, we run the 1-D benchmark problem of \citet{celia1990general} in a Dell Precision 7920 Tower equipped with Intel Xeon Gold 6246R CPU and NVIDIA Quadro RTX 6000 GPU (with 24GB GGDR6 memory). The MP-FVM algorithm is implemented in Python 3.10.5. The total computational time for solving the Celia problem from scratch with $S = 500$ is $181.43$ seconds, in which the neural network training step costs $127.58$ seconds. On the other hand, when using a pre-trained model, the time for neural network training step and the total computational time are reduced by $89.79\%$ and $63.21\%$ down to $13.01$ and $66.76$ seconds, respectively. Meanwhile, the computational time for a direct solver is $43.75$ seconds. While our MP-FVM algorithm still takes more time than the direct solver, it is still an attractive numerical framework as: 1) it gives more accurate solutions; 2) its data-driven nature makes it suitable for seamless integration between physics-based modeling and in situ soil sensing technologies; 3) for large-scale and/or more complex problem settings, the neural network training time will become less significant compared to the actual solution time; and 4) our MP-FVM algorithm consumes less computational time compared to many neural PDE solvers \citep{lu2020,brandstetter}.

\subsection{A 1-D Layered Soil Benchmark Problem}

To investigate the robustness of our MP-FVM algorithm in handling realistic problems, we study the classic Hills' problem \citep{hills1989modeling} that involves the 1-D water infiltration into two layers of very dry soil, each having a depth of 30 cm. The top layer (layer 1) corresponds to Berino loamy fine sand and the bottom layer (layer 2) corresponds to Gledale clay loam. The WRC and HCF follow the Mualem-van Genutchen model. The soil-specific parameters are extracted from \citet{hills1989modeling} and are listed in Table \ref{table_layered_parameters}. This benchmark problem also ignores the sink term.

As pointed out by \citet{berardi20181d}, the dry condition is the most challenging physical case to model from a numerical point of view. The presence of discontinuous interface across the two soil layers presents another complication to this problem. We simulate the problem for up to $7.5$ minutes. For neural network training, we generate a total of 30,500 reference solutions using conventional FVM solver (which implements the static fixed-point iteration scheme of Equation \eqref{eqn_lscheme2} with an optimal $\tau = 0.04$ identified by a trial-and-error procedure).

\begin{table}[ht!]
    \centering
    \begin{tabular}{llllll}
    \toprule
    Soil & $\theta_r$ & $\theta_s$ & $\alpha$ & $n$ & $K_s$\\\midrule
    Berino loamy fine sand & $0.029$ & $0.366$ & $0.028$ & $2.239$ & $541.0$\\
    Gledale clay loam & $0.106$ & $0.469$ & $0.010$ & $1.395$ & $13.10$\\
    \bottomrule
    \end{tabular}
    \caption{Soil-specific parameters and constants used in the layered soil problem of \citet{hills1989modeling}.}
    \label{table_layered_parameters}
\end{table}

\begin{figure}[ht!]
    \centering
    \includegraphics[width=1\textwidth]{fig_compare_layer.jpg}
    \vspace{-2em}
    \caption{Comparison of soil moisture content profile obtained different methods with $\Delta z=1$ cm under (left) MP-FVM, FVM and TMOL at $t=T=3$ sec and $t=T=2.5$ min and (right) MP-FVM, FVM and TMOL at $t=T=7.5$ min. Note that TMOL by \citet{berardi20181d} is not an iterative method. FVM and MP-FVM are implemented for $500$ iterations at every time step.}
    \label{fig_layered_comparison}
\end{figure}

Figure \ref{fig_layered_comparison} illustrates the soil moisture profile at three different times obtained using our MP-FVM algorithm, conventional FVM algorithm, as well as the Transversal Method of Lines (TMOL) solver \citep{berardi20181d} (which is considered the current state-of-the-art algorithm for this problem). All three approaches adopt the same discretized temporal ($\Delta t = 1$ second) and spatial steps ($\Delta z = 1$ cm). We set $\mathrm{RE}_s = 1\times 10^{-5}$ as the common stopping criterion. From Figure \ref{fig_layered_comparison}, we observe that our MP-FVM algorithm is capable of successfully simulating this challenging problem with discontinuities in soil properties at the interface. The soil moisture solutions obtained by our MP-FVM algorithm are also consistent with existing solvers. In fact, compared to the FVM solver, the solutions produced by our MP-FVM algorithm are closer to the state-of-the-art TMOL solutions.

\subsection{A 2-D Benchmark Problem}

In the second example, we study the 2-D Richards equation for an infiltration process in a $1\text{m}\times 1\text{m}$ loam soil field \citep{gkasiorowski2020numerical}. The spatial steps in both horizontal ($\Delta x$) and vertical ($\Delta z$) directions are set to be $0.02$ m, and the time step used for this comparison study is $\Delta t = 10$ seconds. The Mualem-van Genuchten model (see Table \ref{table_hcfwrc}) was used in this case study. The soil-specific parameters, given by \citet{carsel1988developing}, are listed in Table \ref{table_2dparameters}. This problem also ignores the sink term.

\begin{table}[ht!]
    \centering
    \begin{tabular}{llll}
    \toprule
    Property &	Symbol&	Value&	Units\\
    \midrule
    Saturated hydraulic conductivity &$K_s$  & $2.89\times10^{-6}$ & m/s\\
    Saturated water content &$\theta_s$ & $0.43$ & --\\
    Residual water content &$\theta_{r}$ & $0.078$ & --\\
    van Genuchten Constant & $\alpha$ & $3.6$& $\text{m}^{-1}$\\
    van Genuchten Constant & $n$ & $1.56$ & -- \\
    Total time &$T$ & $1.26\times 10^4$ & s\\
    \bottomrule
    \end{tabular}
    \caption{Soil-specific parameters and constants used in 2-D case study. }
    \label{table_2dparameters}
\end{table}

The initial and boundary conditions of this case study are given by:

\begin{subequations}
    \begin{align*}
        &\text{Initial condition: }\psi(x,z,t=0\text{ s})=\left\{
    \begin{aligned}
    0\mathrm{ m}, \quad& x\in [0.46,0.54]\mathrm{ m},\,z=0\mathrm{ m}, \\
    -10\mathrm{ m}, \quad& \text{otherwise.}\\
    \end{aligned}
    \right.\\
    &\text{Boundary condition: }\psi(x\in [0.46, 0.54] \mathrm{ m},z=0,t)=0\mathrm{ m}, \text{no slip conditions for other boundaries.}
    \end{align*}
\end{subequations}

Note that the initial and boundary conditions are symmetric along $x=0.5\mathrm{ m}$. We first obtain $9$ sets of original reference solutions $(\psi, \mu)$, where each $\psi$ or $\mu$ is a $51 \times 51$ array. Here, $\psi$ solutions are obtained from the conventional 2-D FVM solver (which implements the static fixed-point iteration scheme) that uses a spatial step of $0.02$ m under three different fixed-point parameters $\tau = 2,\, 2.22$ and $2.5$ and three total iteration counts $S = 300,\, 400$ and $500$. Then, we apply data augmentation by adding Gaussian noises with $\sigma^2$ values ranging from $0.01$ to $0.05$ to generate a total of $400$ reference solutions (which also contain the original reference solutions). Meanwhile, $\mu$ solutions are obtained from the HYDRUS software \citep{vsimuunek2016recent}. These reference solutions are used to train the encoder-decoder networks for our MP-FVM algorithm. Each neural network contains $3$ hidden layers and $256$ neurons in each layer. ReLU activation function is adopted in each layer, and each neural network is trained by Adam optimizer for 100 epochs. We set the total iteration number to be $S = 500$. The total computational time for our MP-FVM algorithm to run from scratch with $S=500$ is $1473.5$ seconds, whereas the FVM solver takes $876.6$ seconds under the same $S$. 

Meanwhile, we also simulate this 2-D problem using HYDRUS software \citep{vsimuunek2016recent} and compare the pressure head results at $t=T=1.26\times 10^4$ sec with our MP-FVM algorithm and the FVM solver (the fixed-point parameter identified to be $1$ by trial-and-error). From Figure \ref{fig_2dcompare_head}, we can draw two observations. First, the pressure head solution profiles for both FVM and MP-FVM algorithms appear to be symmetric along $x=0.5$ m, whereas HYDRUS 2D shows a clear asymmetric profile. As pointed out earlier, since the initial and boundary conditions are symmetric along $x=0.5$ m, symmetry in the pressure head solutions is expected. This suggests that both FVM and MP-FVM based solvers can capture some degree of underlying physics of the original problem. Second, despite the assymetric behavior in pressure head profile, the size of isolines for the HYDRUS 2D simulation result is more similar to our MP-FVM solution than to the FVM solver solution. This observation is also consistent with the information presented in Figure \ref{fig_2dcompare_1dview}a. In fact, both observations can also be carried over to the soil moisture profile, as shown in Figures \ref{fig_2dcompare_theta} and \ref{fig_2dcompare_1dview}b. Finally, in terms of mass conservation, our MP-FVM algorithm achieves significantly higher $\mathrm{MB}$ value compared to other benchmark solvers (see Table \ref{table_2dmb}).

\begin{figure}[ht!]
    \centering
    \includegraphics[trim = {5cm 1cm 3cm 0}, clip, width=1\textwidth]{fig_2d_compare_psi.jpg}
    \vspace{-2em}
    \caption{Pressure head solution profile obtained from three numerical methods: (left) FVM solver (fixed-point parameter $\tau=1$); (middle) HYDRUS 2D software; (right) our MP-FVM algorithm.}
    \label{fig_2dcompare_head}
\end{figure}

\begin{figure}[ht!]
    \centering
    \includegraphics[trim = {5cm 1cm 3cm 0}, clip, width=\textwidth]{fig_2d_compare_theta.jpg}
    \vspace{-2em}
    \caption{Soil moisture solution profile obtained from three numerical methods: (left) FVM solver (fixed-point parameter $\tau=1$); (middle) HYDRUS 2D software; (right) our MP-FVM algorithm.}
    \label{fig_2dcompare_theta}
\end{figure}

\begin{figure}[ht!]
    \centering
    \includegraphics[width=\textwidth]{fig_2dcompare_1dview.jpg}
    \vspace{-2em}
    \caption{Cross-sectional view ($x=0.5\text{m}$) of: (left) the pressure head profile; (right) soil moisture profile.} \label{fig_2dcompare_1dview}
   \end{figure}

\begin{table}[H]
    \centering
    \begin{tabular}{ll}
    \toprule
    Method &$\mathrm{MB}$ ($\Delta t = 10$ sec) \\ \midrule
    FVM algorithm &$63.12\%$ \\
    HYDRUS 2D simulation & $62.45\%$ \\
    MP-FVM algorithm & \pmb{$71.74\%$}  \\ 
    \bottomrule
    \end{tabular}
    \caption{$\mathrm{MB}$ results of three methods at $x=0.5$ m.}
    \label{table_2dmb}
\end{table}

\subsection{A 3-D Benchmark Problem with Analytical Solutions}

Lastly, we consider a 3-D water infiltration example, in which the analytical solution exists \citep{tracy2006clean}. In this example, $V$ is a 3-D cuboid $[0,a]\times[0,b]\times[0,c]$. The hydraulic conductivity function follows the Gardner's model \citep{gardner1958some} (see Table \ref{table_hcfwrc}). The initial condition is given by:

\begin{equation*}
\psi(x,y,z,t=0)=h_r,
\end{equation*} 
where $h_r$ is a constant. The boundary condition is given by:

\begin{equation*}
    \psi(x,y,z=c,t)=\frac{1}{\alpha}\ln\left[\exp{(\alpha h_r)}+\overline{h}_0\sin\frac{\pi x}{a}\sin\frac{\pi y}{b}\right],
\end{equation*}
where $\overline{h}_0=1-\exp{(\alpha h_r)}$. Ignoring the sink term, the pressure head solution for this problem was derived in \citet{tracy2006clean} as:

\begin{equation}
    \psi=\frac{1}{\alpha}\ln{\bigg\{\exp{(\alpha h_r)}+\overline{h}_0\sin{\frac{\pi x}{a}}\sin{\frac{\pi y}{b}}\exp{\left(\frac{\alpha(c-z)}{2}\right)}\Big[\frac{\sinh{\beta z}}{\sinh{\beta c}}+\frac{2}{z d}\sum_{k=1}^{\infty}(-1)^k\frac{\lambda_k}{\gamma}\sin{(\lambda_k z)}\exp{(-rt)}\Big]\bigg\}},
\end{equation}\label{eqn_3dsolution}
where $d=\frac{\alpha(\theta_s-\theta_r)}{K_s}$, $\lambda_k=\frac{k\pi}{c}$, $\gamma=\frac{\lambda_k^2+\beta^2}{c}$ and $\beta=\sqrt{\frac{\alpha^2}{4}+(\frac{\pi}{a})^2+(\frac{\pi}{b})^2}$.

The infinite series in Equation \eqref{eqn_3dsolution} is convergent by the alternating series test, and we consider the first $1,000$ terms of this series. Note from Equation \eqref{eqn_3dsolution} that the analytical solution depends only on the saturated ($\theta_s$) and residual soil moisture content ($\theta_s$). The Mualem-van Genuchten correlation \citep{Mualem1976,vanGenuchten1980} tabulated in Table \ref{table_hcfwrc} was used for the water retention curve $\theta(\psi)$. The constants and parameters used in this case study are listed in Table \ref{table_3dparameters}.

\begin{table}[ht!]
    \centering
    \begin{tabular}{llll}
    \toprule
    Property &	Symbol&	Value&	Units\\            \midrule
    Saturated hydraulic conductivity &$K_s$  & $1.1$ & m/s\\
    Saturated soil moisture &$\theta_s$ & $0.5$ & --\\
    Residual soil moisture &$\theta_{r}$ & $0$ & --\\
    Parameter in Gardner's model & $\alpha$ & $0.1$ & $\text{m}^{-1}$\\
    Parameter in intial and boundary conditions & $h_r$ & $-15.24$& m\\
    Length of $V$ & $a$ & 2 & m\\
    Width of $V$ & $b$ & 2 & m\\
    Depth of $V$ & $c$ & 2 & m\\
    Total time &$T$ & 86,400 & sec\\
    \bottomrule
    \end{tabular}
    \caption{Soil-specific parameters and constants used in the 3-D case study.}
    \label{table_3dparameters}
\end{table}

Our goal is to compare the accuracy of our MP-FVM algorithm with FVM solvers using this analytical solution as the benchmark. We use our own in-house 3-D FVM solver, which implements the static fixed-point iteration scheme of the FVM-discretized 3-D Richards equation, to obtain $1,734$ original reference solutions using a coarse grid of $\Delta x=\Delta y=\Delta z=0.4$ m under two fixed-point parameters $\tau=1$ and $2$ and five total iteration counts $S = 100,\, 200,\,\dots,\,500$, while excluding any NaN values. Then, data augmentation is applied by introducing Gaussian noise, resulting in a total of $8,820$ data points (which include the original reference solutions) for neural network training. For both FVM and MP-FVM algorithms, we set the tolerance to $1 \times 10^{-9}$, which can be achieved in less than $500$ iterations for each time step. 

\begin{figure}[ht!]
    \centering
    \includegraphics[width=\textwidth]{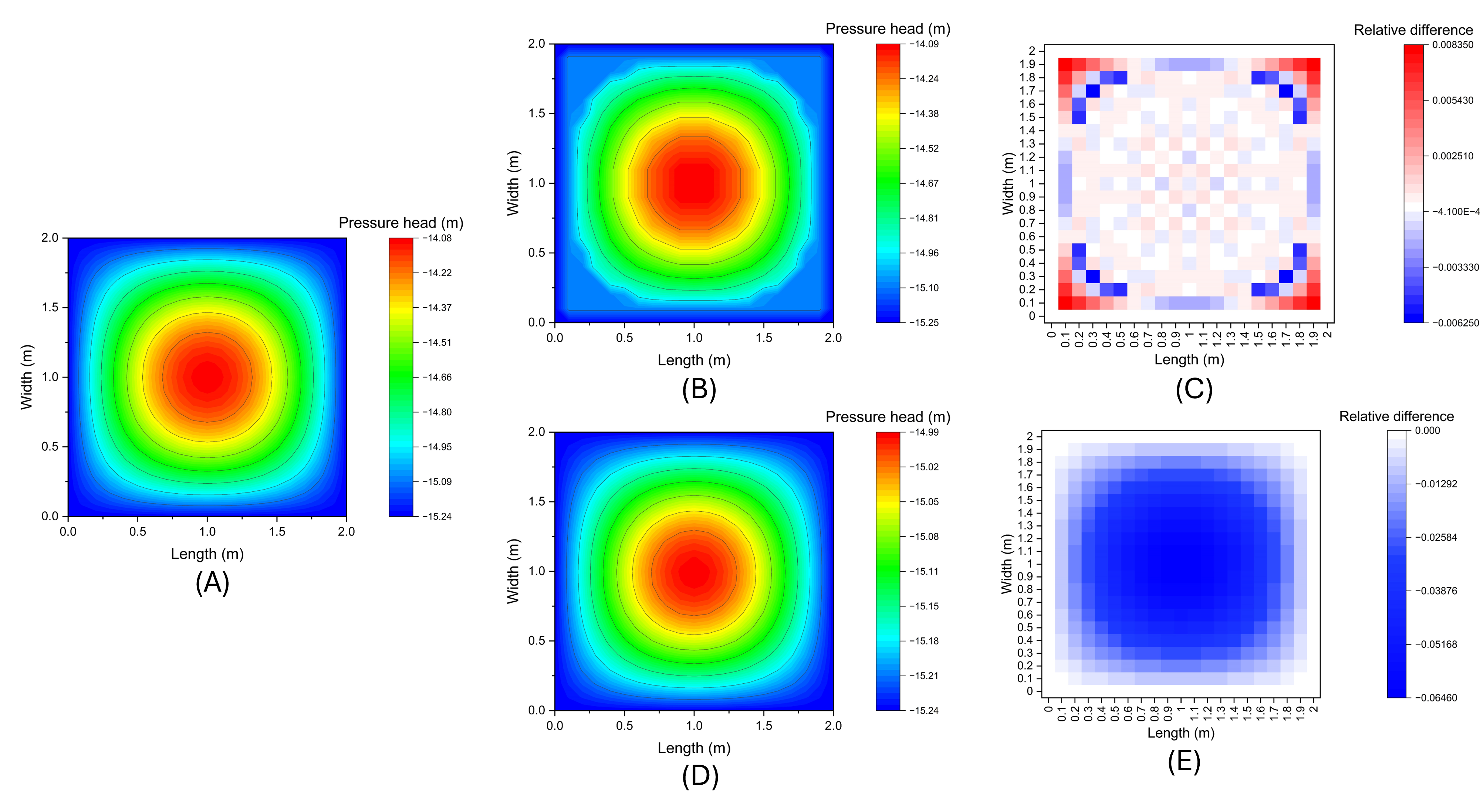}
    \vspace{-2em}
    \caption{Pressure head solution at $z=0.5$ m of different methods: (A) analytical solution, (B) MP-FVM algorithm, (C) the relative difference between analytical and MP-FVM solutions, (D) conventional FVM solver (which implements static fixed-point iteration scheme with an optimal $\tau=2$) and (E) the relative difference between analytical solution and FVM solution.}
    \label{fig_3d0.5m}
\end{figure}

\begin{figure}[ht!]
    \centering
    \includegraphics[width=\textwidth]{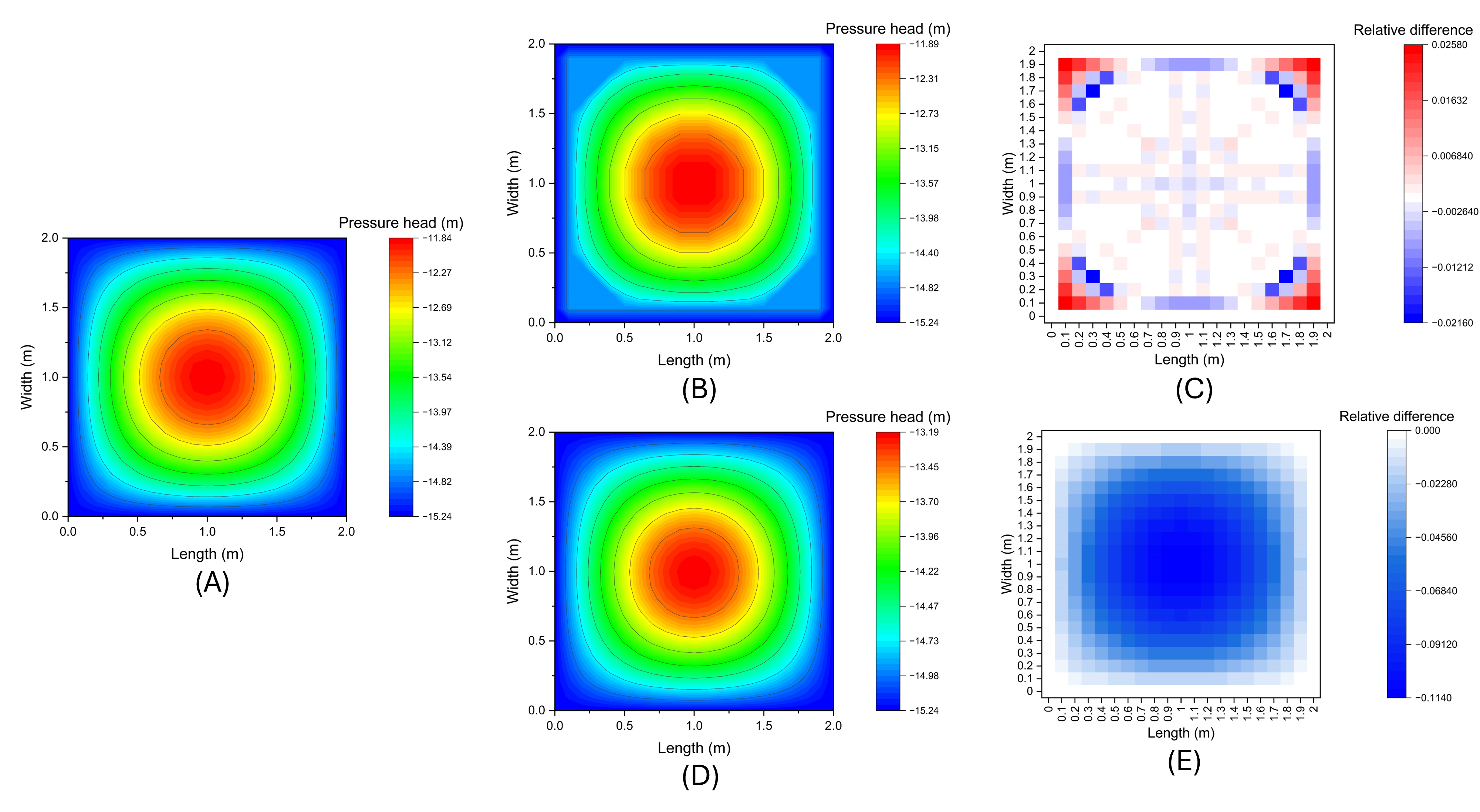}
    \vspace{-2em}
    \caption{Pressure head solution at $z=1.0$ m of different methods: (A) analytical solution, (B) MP-FVM algorithm, (C) the relative difference between analytical and MP-FVM solutions, (D) conventional FVM solver (which implements static fixed-point iteration scheme with an optimal $\tau=2$) and (E) the relative difference between analytical solution and FVM solution.}
    \label{fig_3d1m}
\end{figure}

We examine and compare the pressure head solutions at $z=0.5$ and 1 m, which are shown in Figure \ref{fig_3d0.5m} and \ref{fig_3d1m}, respectively. We quantify the relative difference between the numerical and analytical solutions by $\frac{\psi_{\text{analytical}}-\psi_{\text{numerical}}}{\psi_{\text{analytical}}}$. From the relative difference heat map of Figure \ref{fig_3d0.5m}c,e and \ref{fig_3d1m}c,e, we observe that, first, the magnitude of relative difference of our MP-FVM algorithm is significantly lower than that of the conventional FVM solver. Second, the largest relative difference of our MP-FVM pressure head solution occurs around the four corners of the $x$-$y$ domain, whereas the largest relative difference of FVM solution occurs in the center of the $x$-$y$ domain. Furthermore, in each cell, the relative difference of FVM based pressure head solution is always non-positive, whereas that of MP-FVM based solution can be positive or negative.

Here, we provide some justifications for these observations. First, for conventional FVM solver that embeds the static fixed-point iteration scheme, we observe from Equations \eqref{eqn_lscheme2} that:
\begin{equation*}
    \psi_{\text{analytical}}-\psi_{\text{numerical}}  \propto
    \left\{\sum_{j\in \mathcal{N}_i} \big[K(\psi)\nabla (\psi+z)\big]_{\omega_{i,j}}^{m+1,s}\cdot\mathbf{n}_{\omega_{i,j}}A_{\omega_{i,j}}-\partial_t\theta_i^{m+1}\text{vol}(V_i)\right\},
\end{equation*}
for any $s$, discretized cell $V_i$, and discretized time step $m$. Since the hydraulic conductivity function is positive and symmetric along $x=1$ m and $y=1$ m, and $\nabla \psi\big|_{\omega^+\coloneqq[0,1]\times[0,1]\times z}=-\nabla \psi\big|_{\omega^-\coloneqq[1,2]\times[1,2]\times z}$, we have $\sum_{j\in \mathcal{N}_i} \big[K(\theta(\psi))\nabla (\psi+z)\big]_{\omega_{i,j}}^{m+1,s}\cdot\mathbf{n}_{\omega_{i,j}}A_{\omega_{i,j}} > 0$. Meanwhile, $\partial_t\theta_i^{m+1}(\psi)\text{vol}(V_i)$ is typically small due to the slow dynamics of water infiltration in soil and the fact that $\text{vol}(V_i)$ is small. Thus, we have $\psi_{\text{analytical}}-\psi_{\text{numerical}}>0$ for the FVM solution, which explains why the relative difference is non-positive. On the other hand, for our MP-FVM algorithm, the use of neural networks to approximate $f$ and $f^{-1}$ complicates the behavior (including the sign) of the relative difference.

Regarding the distribution of the magnitude of relative difference in the FVM solver, since hydraulic conductivity function is an increasing function of $\psi$, and $\psi$ is at its maximum at the center of the $x$--$y$ plain, it is expected that $\sum_{j\in \mathcal{N}_{i}} \big[K(\psi)\nabla (\psi+z)\big]_{\omega_{i,j}}^{m+1,s}\cdot\mathbf{n}_{\omega_{i,j}}A_{\omega_{i,j}}$ (hence the relative difference) is maximized at and around the center of the $x$-$y$ plane. However, for MP-FVM based pressure head solution, we suspect that the higher relative difference at the four corners may be attributed to the slight decrease in accuracy of neural networks in approximating $f$ and $f^{-1}$ near the domain boundaries.

Finally, we evaluate the Mean Absolute Error ($\mathrm{MAE}$) by averaging the absolute errors between numerical and analytical pressure head solutions across all cells on two vertical planes, $z=0.5$ m and $z=1$ m. For $z=0.5$ m, $\mathrm{MAE}_{\text{MP-FVM}}$ and $\mathrm{MAE}_{\text{FVM}}$ are calculated to be $0.0146$ and $0.3444$, respectively. For $z=1$ m, $\mathrm{MAE}_{\text{MP-FVM}}$ and $\mathrm{MAE}_{\text{FVM}}$ are $0.0375$ and $0.5653$, respectively. This indicates that the $\mathrm{MAE}$ of the FVM solutions is typically $1$ to $2$ orders of magnitude higher than the MP-FVM solutions, highlighting the accuracy of our MP-FVM algorithm.

\section{A Realistic Case Study} \label{real_case}

Finally, we consider a real-world case study adopted from \citet{orouskhani_impact_2023}, where infiltration, irrigation, and root water extraction take place in circular agricultural field, equipped with a center-pivot irrigation system with a radius of $50$ m, located in Lethbridge, Alberta. Soil moisture sensors are inserted at a depth of 25 cm across 20 different locations in this field to collect soil moisture data every 30 min from June 19 to August 13, 2019. To validate our MP-FVM algorithm in solving real-world 3-D applications, we select one of the 20 locations where the Mualem-van Genuchten WRC and HCF model parameters are identified and given in \citet{orouskhani_impact_2023}. We consider a cylindrical control volume $V$ with a radius of $0.1$ m and a depth of $25$ cm. We discretize $V$ into 6, 40
and 22 nodes in the radial, azimuthal, and axial directions, respectively. The time step size $\Delta t$ is determined using the heuristic formula. Thus, we reformulate Equation \eqref{eqn_DRWlscheme} in cylindrical coordinate system as:
\begin{equation*}
    \mu_i^{m+1,s+1} = \mu_i^{m+1,s}+\tau_i^{m+1,s}\sum_{j\in \mathcal{N}_i}K_{\omega_{i,j}}^{m+1,s}\mathbf{\hat{e}}_j \cdot \mathbf{n}_{\omega_{i,j}}\frac{\mu_j^{m+1,s}-\mu_i^{m+1,s}}{\text{dist}(V_j,V_i)}A_{\omega_{i,j}}+ f^{-1}(J),
\end{equation*}
where $\mathbf{\hat{e}}_j=(1,\frac{1}{r_j^2},1)^T$ and 
\begin{equation}\label{eqn_Jsink}
    \begin{aligned}
        J = & \tau_i^{m+1,s} \sum_{j\in \mathcal{N}_i}K_{\omega_{i,j}}^{m+1,s} \mathbf{\hat{e}}_j \cdot \mathbf{n}_{\omega_{i,j}}\frac{z_j-z_i}{\text{dist}(V_j,V_i)}A_{\omega_{i,j}} - \tau_i^{m+1,s}\frac{\theta_i^{m+1,s}-\theta_i^m}{\Delta t}\text{vol}(V_i) \\ 
        & -\tau_i^{m+1,s}S(\psi_i^{m+1,s})\text{vol}(V_i).
    \end{aligned}
\end{equation}

Here, the sink term in $S$ follows the Feddes model \citep{feddes1978sinkterm}:
\begin{equation}
    S=\sigma(\psi)S_{\text{max}},
\end{equation}
where $S_{\text{max}}$ is the maximum possible root extraction rate and $\sigma$ denotes a dimensionless water stress reduction factor (see \citet{agyeman2021soil} for the detailed formulation).

The boundary conditions are given by:
\begin{equation*}
    \begin{aligned}
        &\frac{\partial \psi(r,\omega,z)}{\partial r}=0 \quad \text{ at }r=0\mathrm{ m},\\
        &\frac{\partial \psi(r,\omega,z)}{\partial r}=0 \quad \text{ at }r=0.1\mathrm{ m},\\
        &\frac{\partial \psi(r,\omega,z)}{\partial z}=0 \quad \text{ at }z=0\text{ cm},\\
        &\frac{\partial \psi(r,\omega,z)}{\partial z}=-1-\frac{u_{\text{irr}}}{K(\psi)} \quad \text{ at }z=25\text{ cm},\\
        &\psi(r,\omega=0,z)=\psi(r,\omega=2\pi,z),
    \end{aligned}
\end{equation*}
where $u_{\text{irr}}$ is the irrigation rate (in m/s). The initial condition is simply:
\begin{equation*}
\psi(x,y,z,t=0)=h_{r},
\end{equation*} 
where $h_r$ is the starting pressure head recording.

Note that the boundary conditions are time dependent due to $u_{\text{irr}}$. This poses a potential computational challenge, as the neural networks typically need to be retrained whenever initial and/or boundary conditions change \citep{mattey2021physics,brecht2023improving}. To overcome this practical challenge, we adopt a new approach of training the two neural networks with 3,000 epochs based on the boundary conditions for June 19, 2019 (no irrigation) when data collection began. Then, the trained weights within these two neural networks serve as the starting point for retraining when a new set of boundary conditions is adopted. This way, only $500$ epochs are sufficient to retrain the neural networks. For each set of boundary conditions, we obtain the training set containing $84,480$ reference solutions. In addition, the dataset provided by \citet{orouskhani_impact_2023}, after performing data augmentation by introducing Gaussian noises, is also included in our training dataset. Each neural netowrk, which has $5$ hidden layers with $256$ neurons in each layer, is trained using SGD optimizer with a learning rate of $0.001$. We set the stopping criterion to be $1 \times 10^{-9}$, which can be achieved well within 500 iterations.

\begin{figure}[ht!]
    \centering
    \includegraphics[trim={0 1cm 0 2cm},clip,width=0.8\textwidth]{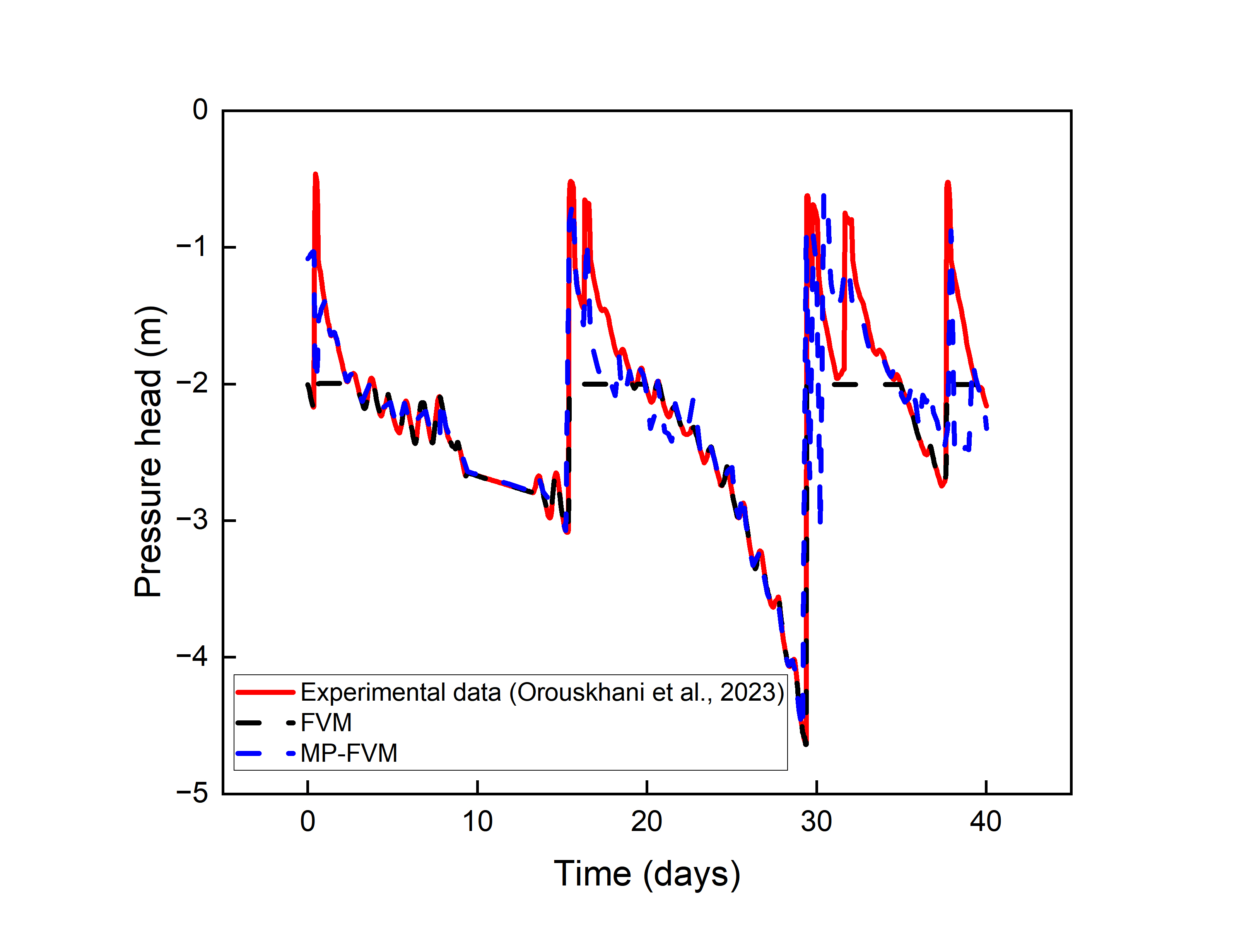}
    \caption{Comparison of pressure head profile at $z = 25$ cm in a selected 0.1-m radius region (averaged for all $6 \times 40 = 240$ cells at $z = 25$ cm) in the field. Note that the standard FVM solver becomes highly inaccurate when the boundary condition changes (15th day, 30th day, etc.). The flattening of true peaks of pressure head solutions represents a nonphysical smoothing of the true solution \citep{miller1998robust}, which we suspect to come from the numerical dispersion and inherent discrete maximum principle (DMP)-type peak clipping behavior observed in standard FVM schemes \citep{njifenjou2025discrete}.}
    \label{fig_real_case}
\end{figure}

For this problem, we simulate the pressure head from 1:00 am on June 19, 2019 to 5:00 pm on July 28, 2019. As mentioned in \citet{orouskhani_impact_2023}, there are two irrigation instances between this time frame, one is on July 4 (the 15th day, $1.81$ mm) and the other is on July 18 (the 30th day, $1.58$ mm). Figure \ref{fig_real_case} shows the pressure head solution profile obtained by our MP-FVM algorithms compared to the experimental measurements provided by \citet{orouskhani_impact_2023} over the course of 35 days. We observe that, most of the time, the MP-FVM solutions match with the experimental measurements very well. The only major mismatch between experimental measurements and MP-FVM solutions occurs on the 30th day, which corresponds to the time when the irrigation takes place. We believe that the mismatch is due to our simplifying assumption regarding the irrigation schedule. Due to the limited information we have on the exact irrigation schedule and intensity, we have to assume that the irrigation instances occurred throughout the day. Thus, we simply divide the irrigation amount by $86,400$ seconds to obtain $u_{\text{irr}}$. However, in reality, the irrigation could end in less than 24 hours. With more accurate $u_{\text{irr}}$ model, our MP-FVM algorithm is expected to produce highly accurate solutions that match more closely with experimental measurements at all times. This makes our MP-FVM algorithm an accurate and scalable numerical framework to solve Richards equation over a long period of time.

\section{Conclusion}\label{conclusion}

In this work, we present a novel message passing finite volume method named the MP-FVM algorithm to accurately and efficiently solve $d$-dimensional Richards equation ($d=1,2,3$). Our MP-FVM algorithm is inspired by the encoder-decoder network architecture \citep{ranade2021discretizationnet} and message passing neural PDE solvers \citep{brandstetter}. It adopts an adaptive fixed-point iteration scheme based on the FVM discretization of the Richards equation and, therefore, significantly improves the convergence and stability of the solution process. To account for the numerical errors observed during actual implementation due to computational constraints and realistic simulation settings, we introduce a data-driven approach to first learn the forward and inverse maps between the solutions using two neural networks, followed by integrating the trained neural networks with the numerical scheme via the message passing mechanism to achieve synergistic improvement in solution quality. Furthermore, we also discuss effective ways, such as the ``coarse-to-fine'' approach and Sobolev training \citep{czarnecki2017}, to perform data augmentation to facilitate neural network training using only a small number of low-fidelity reference solutions as training set. Overall, these innovative techniques work seamlessly to improve the convergence and accuracy of our MP-FVM algorithm in solving the Richards equation. Indeed, via several 1-D through 3-D case studies that span across benchmark problems and real-world applications, we demonstrate that, compared to state-of-the-art numerical solvers, our MP-FVM algorithm not only achieves significantly improved accuracy and convergence, but also better preserves the overall mass balance and conservation laws while being computationally efficient to implement. Finally, we remark that our MP-FVM algorithm is expected to be a generalizable computational framework for modeling a wide range of geotechnical applications, including fractional diffusion \citep{gerolymatou2006modelling, zhang2018long,zheng2025wave,escape_25_fractional}, saturated-unsaturated seepage \citep{sun2024simulation}, in-field monitoring of soil suction profiles \citep{venkatesan2021application}, transport in chemically reactive porous media \citep{saeedmonir2024multiscale, feng2024role}, and so on. 

In terms of future work, we would like to address some of the potential limitations of MP-FVM algorithm. First, since the encoder and decoder networks only approximate the true mappings, small but visible discrepancies may be introduced near the boundaries even when the FVM-based fixed-point iteration scheme by itself matches ground truth solutions. To mitigate this, as discussed previously, we plan to experiment with a hybrid switch-solve approach, where we adopt MP-FVM scheme in the interior of the domain, but fall back to the adaptive fixed-point iteration scheme for boundary cells. Another potential limitation of the current MP-FVM is related to the sensitivity of solution quality to the Sobolev regularization parameter in the loss functions of Equations \eqref{eqn_sobolev} and \eqref{eqn_sobolev_2}. One possible solution is to use staged (homotopy) training, in which one can start by pre-training the model with $\lambda=0$, followed by gradually ramping up $\lambda$ to introduce derivative matching without over-smoothing.

\section*{Acknowledgments}
This research was funded by the U.S. National Science Foundation (NSF) under award number 2442806 and the startup fund of College of Engineering, Architecture, and Technology at Oklahoma State University.

\appendix

\section{Proof of Theorem \ref{theorem_convergence2}}\label{appen_c}
\subsection{Lemmas}
To prove Theorem \ref{theorem_convergence2}, we first need to introduce the following preliminary assumptions and results from \citet{fontaine2021} and \citet{Berner2019}.
\begin{enumerate}[leftmargin=*,label = Assumption \arabic*:]
\item The objective function $f$ is $L$-smooth.
\item There exists a Polish probability space \((Z, \mathcal{Z}, \pi^Z)\) and \(\eta \geq 0\) such that one of the following conditions holds:
\begin{enumerate}[leftmargin = -5em, label = (\alph*)]
    \item There exists a function \(H: \mathbb{R}^d \times Z \to \mathbb{R}^d\) such that for any \(x \in \mathbb{R}^d\), $$\int_Z H(x, z) d\pi^Z(z) = \nabla f(x), \qquad \int_Z \|H(x, z) - \nabla f(x)\|_{L^2}^2 d\pi^Z(z) \leq \eta.$$
    \item There exists a function \(\tilde{f} : \mathbb{R}^d \times Z \to \mathbb{R}\) such that for all \(z \in Z\), \(\tilde{f}(\cdot, z) \in C^1(\mathbb{R}^d, \mathbb{R})\) is \(L\)-smooth. Furthermore, there exists \(x^* \in \mathbb{R}^d\) such that, for any \(x \in \mathbb{R}^d\), $$\int_Z \tilde{f}(x, z) d\pi^Z(z) = f(x), \quad \int_Z \nabla \tilde{f}(x, z) d\pi^Z(z) = \nabla f(x), \quad \int_Z \|\nabla \tilde{f}(x^*, z)\|_{L^2}^2 d\pi^Z(z) \leq \eta.$$ In this case, we define \(H = \nabla \tilde{f}\).
\end{enumerate}
\item There exists \(M \geq 0\) such that for any \(x, y \in \mathbb{R}^d\), \[\|\Sigma(x)^{1/2} - \Sigma(y)^{1/2}\|_{L^2} \leq M \|x - y\|_{L^2}.\]
\item One of the following conditions holds:
\begin{enumerate}[leftmargin = -5em, label = (\alph*)]
    \item For Assumption 2(a): \( f \) is convex, \textit{i.e.}, for any \(x, y \in \mathbb{R}^d\),
    \[\langle \nabla f(x) - \nabla f(y), x - y \rangle \geq 0,\] and there exists a minimizer \(x^* \in \arg \min_{x \in \mathbb{R}^d} f\).
    \item For Assumption 2(b): For all \(z \in Z\), \(\tilde{f}(\cdot, z)\) is convex, and there exists a minimizer \(x^* \in \arg \min_{x \in \mathbb{R}^d} f\).
\end{enumerate}
\end{enumerate}

Under Assumptions 1 and 2, we introduce the sequence $\{X_n\}_{n\in \mathbb{N}}$ starting from $X_0 \in \mathbb{R}^d$ corresponding to SGD with non-increasing step sizes for any $n \in \mathbb{N}$ by:
\begin{equation*}
    X_{n+1} = X_n - \gamma (n+1)^{-\alpha} H(X_n, Z_{n+1}),
\end{equation*}
where \(\gamma > 0\), \(\alpha \in [0, 1]\), and \(\{Z_n\}_{n \in \mathbb{N}}\) is a sequence of independent random variables on a probability space \((\Omega, \mathcal{F}, P)\) valued in \((Z, \mathcal{Z})\) such that for any \(n \in \mathbb{N}\), \(Z_n\) is distributed according to \(\pi^Z\). As \citet{fontaine2021} pointed out, the solution of the following SDE is a continuous counterpart of $\{X_n\}_{n \in \mathbb{N}}$: $$\dif \mathbf{X_t} = - (\gamma + t)^{-\alpha} \nabla f(\mathbf{X_t}) \dif t + \gamma (\gamma + t)^{-2\alpha} \Sigma(\mathbf{X_t})^{1/2} \dif \mathbf{B_t},$$ where \(\gamma_\alpha = \gamma^{1/(1-\alpha)}\) and \((B_t)_{t \geq 0}\) is a \(d\)-dimensional Brownian motion.

Given these preliminaries, we now leverage two established results as lemmas:
\begin{lemma}[Theorem 6 of \citet{fontaine2021}] \label{lemma1}
    Let \(\alpha, \gamma \in (0, 1)\), for \(f \in C^2(\mathbb{R}^d, \mathbb{R})\), there exists \(C \geq 0\)  such that for any \(T \geq 1\), $$\mathbb{E}[f(\mathbf{X_T})] - \min_{x \in \mathbb{R}^d} f \leq C \frac{(1 + \log(T))^2}{T^{\alpha(1-\alpha)}}.$$
\end{lemma}

\begin{lemma}[Equation 35 of \citet{Berner2019}] \label{lemma2}
    Suppose $\tilde{f}$ with an at most polynomially growing derivative is the ``true'' function learned by the neural network. Let $\kappa>0$ be the polynomial growth rate, there exists $D\geq 0$ such that
    \[\|\tilde{f}(x)-\tilde{f}(y)\|_{L^2}\leq D\left(1+\|x\|_{L^2}^{\kappa+2}+\|y\|_{L^2}^{\kappa+2}\right)\|x-y\|_{L^2}\] holds.
\end{lemma}

\subsection{Proof of Theorem \ref{theorem_convergence2}}
With Lemmas \ref{lemma1} and \ref{lemma2}, we are now ready to give the proof of Theorem \ref{theorem_convergence2} which accounts for the convergence of SGD:
\begin{proof}
 To start, we have: 
 \begin{align*}
 \|\mu^{m+1,s+1}-\mu^{m+1,s}\|_{L^2}&\leq \mathbb{E}\left[\|\hat{f}_{\mathrm{NN}}(\psi^{m+1,s+1}, \mathbf{X_{T}})-\hat{f}_{\mathrm{NN}}(\psi^{m+1,s}, \mathbf{X_{T}})\|_{L^2}\right]\\
 &\leq \mathbb{E}\left[\|\hat{f}_{\mathrm{NN}}(\psi^{m+1,s+1}, \mathbf{X_{T}})-\tilde{f}(\psi^{m+1,s+1})\|_{L^2}\right] + \|\tilde{f}(\psi^{m+1,s+1})-\tilde{f}(\psi^{m+1,s})\|_{L^2}\\
 &+ \mathbb{E}\left[\|\tilde{f}(\psi^{m+1,s})-\hat{f}_{\mathrm{NN}}(\psi^{m+1,s}, \mathbf{X_{T}})\|_{L^2}\right],
 \end{align*}
where $\mathbf{X_{T}}$ is the weights of $\hat{f}_{\mathrm{NN}}$ optimized by SGD optimizer, whose process is assumed to be well-approximated by the SDE in Lemma \ref{lemma1}, and $\tilde{f}$ is the true function learned by $\hat{f}_{\mathrm{NN}}$. 

To bound the first and third terms, we define the objective function for the SGD process for a given input $\psi$ as $f_{\psi}(x) = \|\hat{f}_{\mathrm{NN}}(\psi, x) - \tilde{f}(\psi)\|_{L^2}$. We assume this function satisfies the conditions for Lemma \ref{lemma1}. The terms we seek to bound are then precisely of the form $\mathbb{E}[f_{\psi}(\mathbf{X_T})]$. From Lemma \ref{lemma1}, we have:
$$\mathbb{E}[f_{\psi}(\mathbf{X_T})] \leq C \frac{(1 + \log(T))^2}{T^{\alpha(1-\alpha)}} + \min_{x \in \mathbb{R}^d} f_{\psi}(x).$$

We assume the network is a good approximator, such that for a given $\varepsilon > 0$ and for any relevant $\psi$, the minimum error satisfies $\min_{x \in \mathbb{R}^d} f_{\psi}(x) \leq \frac{\varepsilon}{6}$. When the network is trained for a sufficiently large $T$, we can ensure $C \frac{(1 + \log(T))^2}{T^{\alpha(1-\alpha)}} \leq \frac{\varepsilon}{6}$.
Thus, for both the first and third terms, which correspond to $\psi = \psi^{m+1,s+1}$ and $\psi = \psi^{m+1,s}$, we have the bound:
$$\mathbb{E}\left[\|\hat{f}_{\mathrm{NN}}(\psi, \mathbf{X_{T}})-\tilde{f}(\psi)\|_{L^2}\right] \leq \frac{\varepsilon}{6} + \frac{\varepsilon}{6} = \frac{\varepsilon}{3}.$$

Next, for the term $\|\tilde{f}(\psi^{m+1,s+1})-\tilde{f}(\psi^{m+1,s})\|_{L^2}$, it can be bounded using Lemma \ref{lemma2} and the result $\|\psi^{m+1,s+1}-\psi^{m+1,s}\|_{L^2}\leq \frac{\varepsilon}{3D\left(1+\|\psi^{m+1,s+1}\|_{L^2}^{\kappa+2}+\|\psi^{m+1,s}\|_{L^2}^{\kappa+2}\right)}$ obtained from Theorem \ref{theorem_convergence}:
\begin{align*}\nonumber
 \|\tilde{f}(\psi^{m+1,s+1})-\tilde{f}(\psi^{m+1,s})\|_{L^2}&\leq D\left(1+\|\psi^{m+1,s+1}\|_{L^2}^{\kappa+2}+\|\psi^{m+1,s}\|_{L^2}^{\kappa+2}\right)\|\psi^{m+1,s+1}-\psi^{m+1,s}\|_{L^2}\\
 & \leq \frac{\varepsilon}{3}.
\end{align*}

Therefore, it follows that $$\|\mu^{m+1,s+1}-\mu^{m+1,s}\|_{L^2}\leq \frac{\varepsilon}{3}+\frac{\varepsilon}{3}+\frac{\varepsilon}{3}=\varepsilon, $$ which completes the proof.
\end{proof}


\begin{thebibliography}{82}
\expandafter\ifx\csname natexlab\endcsname\relax\def\natexlab#1{#1}\fi
\providecommand{\url}[1]{\texttt{#1}}
\providecommand{\href}[2]{#2}
\providecommand{\path}[1]{#1}
\providecommand{\DOIprefix}{doi:}
\providecommand{\ArXivprefix}{arXiv:}
\providecommand{\URLprefix}{URL: }
\providecommand{\Pubmedprefix}{pmid:}
\providecommand{\doi}[1]{\href{http://dx.doi.org/#1}{\path{#1}}}
\providecommand{\Pubmed}[1]{\href{pmid:#1}{\path{#1}}}
\providecommand{\bibinfo}[2]{#2}
\ifx\xfnm\relax \def\xfnm[#1]{\unskip,\space#1}\fi
\bibitem[{Abdellatif et~al.(2018)Abdellatif, Bernardi, Touihri and Yakoubi}]{abdellatif2016priori}
\bibinfo{author}{Abdellatif, N.}, \bibinfo{author}{Bernardi, C.}, \bibinfo{author}{Touihri, M.}, \bibinfo{author}{Yakoubi, D.}, \bibinfo{year}{2018}.
\newblock \bibinfo{title}{A priori error analysis of the implicit {E}uler, spectral discretization of a nonlinear equation for a flow in a partially saturated porous media}.
\newblock \bibinfo{journal}{Advances in Pure and Applied Mathematics} \bibinfo{volume}{9}, \bibinfo{pages}{1--27}.
\bibitem[{Agyeman et~al.(2020)Agyeman, Bo, Sahoo, Yin, Liu and Shah}]{agyeman2021soil}
\bibinfo{author}{Agyeman, B.T.}, \bibinfo{author}{Bo, S.}, \bibinfo{author}{Sahoo, S.R.}, \bibinfo{author}{Yin, X.}, \bibinfo{author}{Liu, J.}, \bibinfo{author}{Shah, S.L.}, \bibinfo{year}{2020}.
\newblock \bibinfo{title}{Soil moisture map construction using microwave remote sensors and sequential data assimilation}.
\newblock \URLprefix \url{https://arxiv.org/abs/2010.07037}, \href{http://arxiv.org/abs/2010.07037}{{\tt arXiv:2010.07037}}.
\bibitem[{Amrein(2019)}]{amrein2019adaptive}
\bibinfo{author}{Amrein, M.}, \bibinfo{year}{2019}.
\newblock \bibinfo{title}{Adaptive fixed point iterations for semilinear elliptic partial differential equations}.
\newblock \bibinfo{journal}{Calcolo} \bibinfo{volume}{56}, \bibinfo{pages}{30}.
\bibitem[{Bandai and Ghezzehei(2021)}]{pinn}
\bibinfo{author}{Bandai, T.}, \bibinfo{author}{Ghezzehei, T.A.}, \bibinfo{year}{2021}.
\newblock \bibinfo{title}{Physics-informed neural networks with monotonicity constraints for {R}ichardson-{R}ichards equation: Estimation of constitutive relationships and soil water flux density from volumetric water content measurements}.
\newblock \bibinfo{journal}{Water Resources Research} \bibinfo{volume}{57}, \bibinfo{pages}{e2020WR027642}.
\bibitem[{Bar-Sinai et~al.(2019)Bar-Sinai, Hoyer, Hickey and Brenner}]{bar2019learning}
\bibinfo{author}{Bar-Sinai, Y.}, \bibinfo{author}{Hoyer, S.}, \bibinfo{author}{Hickey, J.}, \bibinfo{author}{Brenner, M.P.}, \bibinfo{year}{2019}.
\newblock \bibinfo{title}{Learning data-driven discretizations for partial differential equations}.
\newblock \bibinfo{journal}{Proceedings of the National Academy of Sciences} \bibinfo{volume}{116}, \bibinfo{pages}{15344--15349}.
\bibitem[{Bassetto et~al.(2022)Bassetto, Canc{\`e}s, Ench{\'e}ry and Tran}]{bassetto2022several}
\bibinfo{author}{Bassetto, S.}, \bibinfo{author}{Canc{\`e}s, C.}, \bibinfo{author}{Ench{\'e}ry, G.}, \bibinfo{author}{Tran, Q.H.}, \bibinfo{year}{2022}.
\newblock \bibinfo{title}{On several numerical strategies to solve {Richards’} equation in heterogeneous media with finite volumes}.
\newblock \bibinfo{journal}{Computational Geosciences} \bibinfo{volume}{26}, \bibinfo{pages}{1297--1322}.
\bibitem[{Belfort et~al.(2013)Belfort, Younes, Fahs and Lehmann}]{oscillation}
\bibinfo{author}{Belfort, B.}, \bibinfo{author}{Younes, A.}, \bibinfo{author}{Fahs, M.}, \bibinfo{author}{Lehmann, F.}, \bibinfo{year}{2013}.
\newblock \bibinfo{title}{On equivalent hydraulic conductivity for oscillation-free solutions of {R}ichard's equation}.
\newblock \bibinfo{journal}{Journal of Hydrology} \bibinfo{volume}{505}, \bibinfo{pages}{202--217}.
\bibitem[{Berardi et~al.(2018)Berardi, Difonzo, Vurro and Lopez}]{berardi20181d}
\bibinfo{author}{Berardi, M.}, \bibinfo{author}{Difonzo, F.}, \bibinfo{author}{Vurro, M.}, \bibinfo{author}{Lopez, L.}, \bibinfo{year}{2018}.
\newblock \bibinfo{title}{The 1{D} {Richards}' equation in two layered soils: a filippov approach to treat discontinuities}.
\newblock \bibinfo{journal}{Advances in Water Resources} \bibinfo{volume}{115}, \bibinfo{pages}{264--272}.
\bibitem[{Bergamaschi and Putti(1999)}]{bergamaschi1999mixed}
\bibinfo{author}{Bergamaschi, L.}, \bibinfo{author}{Putti, M.}, \bibinfo{year}{1999}.
\newblock \bibinfo{title}{Mixed finite elements and {Newton}-type linearizations for the solution of {Richards'} equation}.
\newblock \bibinfo{journal}{International Journal for Numerical Methods in Engineering} \bibinfo{volume}{45}, \bibinfo{pages}{1025--1046}.
\bibitem[{Berner et~al.(2019)Berner, Elbrächter, Grohs and Jentzen}]{Berner2019}
\bibinfo{author}{Berner, J.}, \bibinfo{author}{Elbrächter, D.}, \bibinfo{author}{Grohs, P.}, \bibinfo{author}{Jentzen, A.}, \bibinfo{year}{2019}.
\newblock \bibinfo{title}{Towards a regularity theory for {ReLU} networks – chain rule and global error estimates}, in: \bibinfo{booktitle}{2019 13th International conference on Sampling Theory and Applications (SampTA)}, pp. \bibinfo{pages}{1--5}.
\newblock \DOIprefix\doi{10.1109/SampTA45681.2019.9031005}.
\bibitem[{Brandstetter et~al.(2022)Brandstetter, Worrall and Welling}]{brandstetter}
\bibinfo{author}{Brandstetter, J.}, \bibinfo{author}{Worrall, D.E.}, \bibinfo{author}{Welling, M.}, \bibinfo{year}{2022}.
\newblock \bibinfo{title}{Message passing neural {PDE} solvers}.
\newblock \bibinfo{journal}{CoRR} \bibinfo{volume}{abs/2202.03376}.
\newblock \URLprefix \url{https://arxiv.org/abs/2202.03376}, \href{http://arxiv.org/abs/2202.03376}{{\tt arXiv:2202.03376}}.
\bibitem[{Brecht et~al.(2023)Brecht, Bakels, Bihlo and Stohl}]{brecht2023improving}
\bibinfo{author}{Brecht, R.}, \bibinfo{author}{Bakels, L.}, \bibinfo{author}{Bihlo, A.}, \bibinfo{author}{Stohl, A.}, \bibinfo{year}{2023}.
\newblock \bibinfo{title}{Improving trajectory calculations by flexpart 10.4+ using single-image super-resolution}.
\newblock \bibinfo{journal}{Geoscientific Model Development} \bibinfo{volume}{16}, \bibinfo{pages}{2181--2192}.
\bibitem[{Caputo and Stepanyants(2008)}]{caputo_front_2008}
\bibinfo{author}{Caputo, J.G.}, \bibinfo{author}{Stepanyants, Y.A.}, \bibinfo{year}{2008}.
\newblock \bibinfo{title}{Front solutions of {Richards}' equation}.
\newblock \bibinfo{journal}{Transport in Porous Media} \bibinfo{volume}{74}, \bibinfo{pages}{1--20}.
\bibitem[{Carsel and Parrish(1988)}]{carsel1988developing}
\bibinfo{author}{Carsel, R.F.}, \bibinfo{author}{Parrish, R.S.}, \bibinfo{year}{1988}.
\newblock \bibinfo{title}{Developing joint probability distributions of soil water retention characteristics}.
\newblock \bibinfo{journal}{Water Resources Research} \bibinfo{volume}{24}, \bibinfo{pages}{755--769}.
\bibitem[{Casulli and Zanolli(2010)}]{casulli2010nested}
\bibinfo{author}{Casulli, V.}, \bibinfo{author}{Zanolli, P.}, \bibinfo{year}{2010}.
\newblock \bibinfo{title}{A nested {Newton-type} algorithm for finite volume methods solving {Richards'} equation in mixed form}.
\newblock \bibinfo{journal}{SIAM Journal on Scientific Computing} \bibinfo{volume}{32}, \bibinfo{pages}{2255--2273}.
\bibitem[{Caviedes-Voulli{\`e}me et~al.(2013)Caviedes-Voulli{\`e}me, Garc{\i}, Murillo et~al.}]{caviedes2013verification}
\bibinfo{author}{Caviedes-Voulli{\`e}me, D.}, \bibinfo{author}{Garc{\i}, P.}, \bibinfo{author}{Murillo, J.}, et~al., \bibinfo{year}{2013}.
\newblock \bibinfo{title}{Verification, conservation, stability and efficiency of a finite volume method for the {1D Richards} equation}.
\newblock \bibinfo{journal}{Journal of Hydrology} \bibinfo{volume}{480}, \bibinfo{pages}{69--84}.
\bibitem[{Celia et~al.(1990)Celia, Bouloutas and Zarba}]{celia1990general}
\bibinfo{author}{Celia, M.A.}, \bibinfo{author}{Bouloutas, E.T.}, \bibinfo{author}{Zarba, R.L.}, \bibinfo{year}{1990}.
\newblock \bibinfo{title}{A general mass-conservative numerical solution for the unsaturated flow equation}.
\newblock \bibinfo{journal}{Water Resources Research} \bibinfo{volume}{26}, \bibinfo{pages}{1483--1496}.
\bibitem[{Celia and Zarba(1988)}]{zarba1988numerical}
\bibinfo{author}{Celia, M.A.}, \bibinfo{author}{Zarba, R.}, \bibinfo{year}{1988}.
\newblock \bibinfo{title}{A comparative study of numerical solutions for unsaturated flow}, in: \bibinfo{editor}{Atluri, S.N.}, \bibinfo{editor}{Yagawa, G.} (Eds.), \bibinfo{booktitle}{Computational Mechanics '88}, \bibinfo{publisher}{Springer Berlin Heidelberg}, \bibinfo{address}{Berlin, Heidelberg}. pp. \bibinfo{pages}{1659--1662}.
\bibitem[{Ch{\'a}vez-Negrete et~al.(2024)Ch{\'a}vez-Negrete, Dom{\'\i}nguez-Mota and Rom{\'a}n-Guti{\'e}rrez}]{chavez2024interface}
\bibinfo{author}{Ch{\'a}vez-Negrete, C.}, \bibinfo{author}{Dom{\'\i}nguez-Mota, F.}, \bibinfo{author}{Rom{\'a}n-Guti{\'e}rrez, R.}, \bibinfo{year}{2024}.
\newblock \bibinfo{title}{Interface formulation for generalized finite difference method for solving groundwater flow}.
\newblock \bibinfo{journal}{Computers and Geotechnics} \bibinfo{volume}{166}, \bibinfo{pages}{105990}.
\bibitem[{Chen et~al.(2023)Chen, Xu, Wang and Li}]{chen2023modeling}
\bibinfo{author}{Chen, Y.}, \bibinfo{author}{Xu, Y.}, \bibinfo{author}{Wang, L.}, \bibinfo{author}{Li, T.}, \bibinfo{year}{2023}.
\newblock \bibinfo{title}{Modeling water flow in unsaturated soils through physics-informed neural network with principled loss function}.
\newblock \bibinfo{journal}{Computers and Geotechnics} \bibinfo{volume}{161}, \bibinfo{pages}{105546}.
\bibitem[{Czarnecki et~al.(2017)Czarnecki, Osindero, Jaderberg, Swirszcz and Pascanu}]{czarnecki2017}
\bibinfo{author}{Czarnecki, W.M.}, \bibinfo{author}{Osindero, S.}, \bibinfo{author}{Jaderberg, M.}, \bibinfo{author}{Swirszcz, G.}, \bibinfo{author}{Pascanu, R.}, \bibinfo{year}{2017}.
\newblock \bibinfo{title}{Sobolev training for neural networks}, in: \bibinfo{booktitle}{Proceedings of the 31st International Conference on Neural Information Processing Systems}, \bibinfo{address}{Red Hook, NY, USA}. p. \bibinfo{pages}{4281–4290}.
\bibitem[{Da~Silva and Adeodato(2011)}]{da2011pca}
\bibinfo{author}{Da~Silva, I.B.V.}, \bibinfo{author}{Adeodato, P.J.}, \bibinfo{year}{2011}.
\newblock \bibinfo{title}{{PCA} and {Gaussian} noise in {MLP} neural network training improve generalization in problems with small and unbalanced data sets}, in: \bibinfo{booktitle}{The 2011 International Joint Conference on Neural Networks}, \bibinfo{organization}{IEEE}. pp. \bibinfo{pages}{2664--2669}.
\bibitem[{Day and Luthin(1956)}]{day1956numerical}
\bibinfo{author}{Day, P.R.}, \bibinfo{author}{Luthin, J.N.}, \bibinfo{year}{1956}.
\newblock \bibinfo{title}{A numerical solution of the differential equation of flow for a vertical drainage problem}.
\newblock \bibinfo{journal}{Soil Science Society of America Journal} \bibinfo{volume}{20}, \bibinfo{pages}{443--447}.
\bibitem[{De~Moura and Kubrusly(2013)}]{de2013courant}
\bibinfo{author}{De~Moura, C.A.}, \bibinfo{author}{Kubrusly, C.S.}, \bibinfo{year}{2013}.
\newblock \bibinfo{title}{{The Courant--Friedrichs--Lewy (CFL) condition}}.
\newblock \bibinfo{journal}{AMC} \bibinfo{volume}{10}, \bibinfo{pages}{45--90}.
\bibitem[{Edelsbrunner and Morozov(2013)}]{edelsbrunner2013persistent}
\bibinfo{author}{Edelsbrunner, H.}, \bibinfo{author}{Morozov, D.}, \bibinfo{year}{2013}.
\newblock \bibinfo{title}{Persistent homology: theory and practice}.
\newblock \bibinfo{publisher}{eScholarship, University of California}.
\bibitem[{Evans(2010)}]{evans2010partial}
\bibinfo{author}{Evans, L.C.}, \bibinfo{year}{2010}.
\newblock \bibinfo{title}{Partial differential equations}. volume~\bibinfo{volume}{19}.
\newblock \bibinfo{publisher}{American Mathematical Soc.}
\bibitem[{Eymard et~al.(1999)Eymard, Gutnic and Hilhorst}]{eymard1999finite}
\bibinfo{author}{Eymard, R.}, \bibinfo{author}{Gutnic, M.}, \bibinfo{author}{Hilhorst, D.}, \bibinfo{year}{1999}.
\newblock \bibinfo{title}{The finite volume method for {Richards} equation}.
\newblock \bibinfo{journal}{Computational Geosciences} \bibinfo{volume}{3}, \bibinfo{pages}{259--294}.
\bibitem[{Farthing and Ogden(2017a)}]{farthing2017numerical}
\bibinfo{author}{Farthing, M.W.}, \bibinfo{author}{Ogden, F.L.}, \bibinfo{year}{2017}a.
\newblock \bibinfo{title}{Numerical solution of {Richards'} equation: A review of advances and challenges}.
\newblock \bibinfo{journal}{Soil Science Society of America Journal} \bibinfo{volume}{81}, \bibinfo{pages}{1257--1269}.
\bibitem[{Farthing and Ogden(2017b)}]{richardsreview}
\bibinfo{author}{Farthing, M.W.}, \bibinfo{author}{Ogden, F.L.}, \bibinfo{year}{2017}b.
\newblock \bibinfo{title}{Numerical solution of {Richards'} equation: A review of advances and challenges}.
\newblock \bibinfo{journal}{Soil Science Society of America journal} \bibinfo{volume}{81}, \bibinfo{pages}{1257--1269}.
\bibitem[{Feddes and Zaradny(1978)}]{feddes1978sinkterm}
\bibinfo{author}{Feddes, R.}, \bibinfo{author}{Zaradny, H.}, \bibinfo{year}{1978}.
\newblock \bibinfo{title}{Model for simulating soil-water content considering evapotranspiration --- comments}.
\newblock \bibinfo{journal}{Journal of Hydrology} \bibinfo{volume}{37}, \bibinfo{pages}{393--397}.
\bibitem[{Feng et~al.(2024)Feng, Wang, Wang, Chen and Chen}]{feng2024role}
\bibinfo{author}{Feng, S.J.}, \bibinfo{author}{Wang, H.Y.}, \bibinfo{author}{Wang, W.T.}, \bibinfo{author}{Chen, H.}, \bibinfo{author}{Chen, H.X.}, \bibinfo{year}{2024}.
\newblock \bibinfo{title}{Role of physical and chemical nonequilibriums in isco remediation of contaminated soils}.
\newblock \bibinfo{journal}{Computers and Geotechnics} \bibinfo{volume}{176}, \bibinfo{pages}{106725}.
\bibitem[{Fontaine et~al.(2021)Fontaine, Bortoli and Durmus}]{fontaine2021}
\bibinfo{author}{Fontaine, X.}, \bibinfo{author}{Bortoli, V.D.}, \bibinfo{author}{Durmus, A.}, \bibinfo{year}{2021}.
\newblock \bibinfo{title}{Convergence rates and approximation results for {SGD} and its continuous-time counterpart}, in: \bibinfo{editor}{Belkin, M.}, \bibinfo{editor}{Kpotufe, S.} (Eds.), \bibinfo{booktitle}{Proceedings of Thirty Fourth Conference on Learning Theory}, \bibinfo{publisher}{PMLR}. pp. \bibinfo{pages}{1965--2058}.
\bibitem[{Gardner(1958)}]{gardner1958some}
\bibinfo{author}{Gardner, W.}, \bibinfo{year}{1958}.
\newblock \bibinfo{title}{Some steady-state solutions of the unsaturated moisture flow equation with application to evaporation from a water table}.
\newblock \bibinfo{journal}{Soil Science} \bibinfo{volume}{85}, \bibinfo{pages}{228--232}.
\bibitem[{G{\k{a}}siorowski and Kolerski(2020)}]{gkasiorowski2020numerical}
\bibinfo{author}{G{\k{a}}siorowski, D.}, \bibinfo{author}{Kolerski, T.}, \bibinfo{year}{2020}.
\newblock \bibinfo{title}{Numerical solution of the two-dimensional {Richards} equation using alternate splitting methods for dimensional decomposition}.
\newblock \bibinfo{journal}{Water} \bibinfo{volume}{12}, \bibinfo{pages}{1780}.
\bibitem[{Gerolymatou et~al.(2006)Gerolymatou, Vardoulakis and Hilfer}]{gerolymatou2006modelling}
\bibinfo{author}{Gerolymatou, E.}, \bibinfo{author}{Vardoulakis, I.}, \bibinfo{author}{Hilfer, R.}, \bibinfo{year}{2006}.
\newblock \bibinfo{title}{Modelling infiltration by means of a nonlinear fractional diffusion model}.
\newblock \bibinfo{journal}{Journal of Physics D: Applied Physics} \bibinfo{volume}{39}, \bibinfo{pages}{4104}.
\bibitem[{Haghighat et~al.(2023)Haghighat, Sattari and Wuttke}]{haghighat2023finite}
\bibinfo{author}{Haghighat, N.}, \bibinfo{author}{Sattari, A.S.}, \bibinfo{author}{Wuttke, F.}, \bibinfo{year}{2023}.
\newblock \bibinfo{title}{Finite discrete element modeling of desiccation fracturing in partially saturated porous medium}.
\newblock \bibinfo{journal}{Computers and Geotechnics} \bibinfo{volume}{164}, \bibinfo{pages}{105761}.
\bibitem[{Haverkamp et~al.(1977)Haverkamp, Vauclin, Touma, Wierenga and Vachaud}]{haverkamp1977comparison}
\bibinfo{author}{Haverkamp, R.}, \bibinfo{author}{Vauclin, M.}, \bibinfo{author}{Touma, J.}, \bibinfo{author}{Wierenga, P.}, \bibinfo{author}{Vachaud, G.}, \bibinfo{year}{1977}.
\newblock \bibinfo{title}{A comparison of numerical simulation models for one-dimensional infiltration}.
\newblock \bibinfo{journal}{Soil Science Society of America Journal} \bibinfo{volume}{41}, \bibinfo{pages}{285--294}.
\bibitem[{Hills et~al.(1989)Hills, Porro, Hudson and Wierenga}]{hills1989modeling}
\bibinfo{author}{Hills, R.G.}, \bibinfo{author}{Porro, I.}, \bibinfo{author}{Hudson, D.B.}, \bibinfo{author}{Wierenga, P.J.}, \bibinfo{year}{1989}.
\newblock \bibinfo{title}{Modeling one-dimensional infiltration into very dry soils: 1. model development and evaluation}.
\newblock \bibinfo{journal}{Water Resources Research} \bibinfo{volume}{25}, \bibinfo{pages}{1259--1269}.
\bibitem[{Hornik(1991)}]{hornik1991approximation}
\bibinfo{author}{Hornik, K.}, \bibinfo{year}{1991}.
\newblock \bibinfo{title}{Approximation capabilities of multilayer feedforward networks}.
\newblock \bibinfo{journal}{Neural Networks} \bibinfo{volume}{4}, \bibinfo{pages}{251--257}.
\bibitem[{Hyman and Shashkov(1997)}]{hyman1997natural}
\bibinfo{author}{Hyman, J.M.}, \bibinfo{author}{Shashkov, M.}, \bibinfo{year}{1997}.
\newblock \bibinfo{title}{Natural discretizations for the divergence, gradient, and curl on logically rectangular grids}.
\newblock \bibinfo{journal}{Computers \& Mathematics with Applications} \bibinfo{volume}{33}, \bibinfo{pages}{81--104}.
\bibitem[{Ireson et~al.(2023)Ireson, Spiteri, Clark and Mathias}]{ireson2023simple}
\bibinfo{author}{Ireson, A.M.}, \bibinfo{author}{Spiteri, R.J.}, \bibinfo{author}{Clark, M.P.}, \bibinfo{author}{Mathias, S.A.}, \bibinfo{year}{2023}.
\newblock \bibinfo{title}{A simple, efficient, mass-conservative approach to solving {Richards'} equation {(openRE, v1. 0)}}.
\newblock \bibinfo{journal}{Geoscientific Model Development} \bibinfo{volume}{16}, \bibinfo{pages}{659--677}.
\bibitem[{Lai and Ogden(2015)}]{lai2015mass}
\bibinfo{author}{Lai, W.}, \bibinfo{author}{Ogden, F.L.}, \bibinfo{year}{2015}.
\newblock \bibinfo{title}{A mass-conservative finite volume predictor--corrector solution of the {1D Richards’} equation}.
\newblock \bibinfo{journal}{Journal of Hydrology} \bibinfo{volume}{523}, \bibinfo{pages}{119--127}.
\bibitem[{Lan et~al.(2024)Lan, Su, Zhu, Huang and Zhang}]{lan2024reconstructing}
\bibinfo{author}{Lan, P.}, \bibinfo{author}{Su, J.}, \bibinfo{author}{Zhu, S.}, \bibinfo{author}{Huang, J.}, \bibinfo{author}{Zhang, S.}, \bibinfo{year}{2024}.
\newblock \bibinfo{title}{Reconstructing unsaturated infiltration behavior with sparse data via physics-informed deep learning}.
\newblock \bibinfo{journal}{Computers and Geotechnics} \bibinfo{volume}{168}, \bibinfo{pages}{106162}.
\bibitem[{Lu et~al.(2020)Lu, Kim and Solja\ifmmode \check{c}\else \v{c}\fi{}i\ifmmode~\acute{c}\else \'{c}\fi{}}]{lu2020}
\bibinfo{author}{Lu, P.Y.}, \bibinfo{author}{Kim, S.}, \bibinfo{author}{Solja\ifmmode \check{c}\else \v{c}\fi{}i\ifmmode~\acute{c}\else \'{c}\fi{}, M.}, \bibinfo{year}{2020}.
\newblock \bibinfo{title}{Extracting interpretable physical parameters from spatiotemporal systems using unsupervised learning}.
\newblock \bibinfo{journal}{Phys. Rev. X} \bibinfo{volume}{10}, \bibinfo{pages}{031056}.
\newblock \URLprefix \url{https://link.aps.org/doi/10.1103/PhysRevX.10.031056}.
\bibitem[{Manzini and Ferraris(2004)}]{manzini2004mass}
\bibinfo{author}{Manzini, G.}, \bibinfo{author}{Ferraris, S.}, \bibinfo{year}{2004}.
\newblock \bibinfo{title}{Mass-conservative finite volume methods on {2-D} unstructured grids for the {Richards’} equation}.
\newblock \bibinfo{journal}{Advances in Water Resources} \bibinfo{volume}{27}, \bibinfo{pages}{1199--1215}.
\bibitem[{Mattey and Ghosh(2022)}]{mattey2021physics}
\bibinfo{author}{Mattey, R.}, \bibinfo{author}{Ghosh, S.}, \bibinfo{year}{2022}.
\newblock \bibinfo{title}{A novel sequential method to train physics informed neural networks for {Allen Cahn and Cahn Hilliard} equations}.
\newblock \bibinfo{journal}{Computer Methods in Applied Mechanics and Engineering} \bibinfo{volume}{390}, \bibinfo{pages}{114474}.
\newblock \URLprefix \url{http://dx.doi.org/10.1016/j.cma.2021.114474}, \DOIprefix\doi{10.1016/j.cma.2021.114474}.
\bibitem[{Merz and Rybka(2010)}]{merz2010strong}
\bibinfo{author}{Merz, W.}, \bibinfo{author}{Rybka, P.}, \bibinfo{year}{2010}.
\newblock \bibinfo{title}{Strong solutions to the {R}ichards equation in the unsaturated zone}.
\newblock \bibinfo{journal}{Journal of Mathematical Analysis and Applications} \bibinfo{volume}{371}, \bibinfo{pages}{741--749}.
\bibitem[{Miller et~al.(1998)Miller, Williams, Kelley and Tocci}]{miller1998robust}
\bibinfo{author}{Miller, C.T.}, \bibinfo{author}{Williams, G.A.}, \bibinfo{author}{Kelley, C.T.}, \bibinfo{author}{Tocci, M.D.}, \bibinfo{year}{1998}.
\newblock \bibinfo{title}{Robust solution of richards' equation for nonuniform porous media}.
\newblock \bibinfo{journal}{Water Resources Research} \bibinfo{volume}{34}, \bibinfo{pages}{2599--2610}.
\bibitem[{Misiats and Lipnikov(2013)}]{misiats2013second}
\bibinfo{author}{Misiats, O.}, \bibinfo{author}{Lipnikov, K.}, \bibinfo{year}{2013}.
\newblock \bibinfo{title}{Second-order accurate monotone finite volume scheme for {Richards’} equation}.
\newblock \bibinfo{journal}{Journal of Computational Physics} \bibinfo{volume}{239}, \bibinfo{pages}{123--137}.
\bibitem[{Mualem(1976)}]{Mualem1976}
\bibinfo{author}{Mualem, Y.}, \bibinfo{year}{1976}.
\newblock \bibinfo{title}{A new model for predicting the hydraulic conductivity of unsaturated porous media}.
\newblock \bibinfo{journal}{Water Resources Research} \bibinfo{volume}{12}, \bibinfo{pages}{513--522}.
\bibitem[{Ng et~al.(2025)Ng, Zhou and Zhang}]{ng2025novel}
\bibinfo{author}{Ng, C.W.}, \bibinfo{author}{Zhou, Q.}, \bibinfo{author}{Zhang, Q.}, \bibinfo{year}{2025}.
\newblock \bibinfo{title}{A novel surrogate model for hydro-mechanical coupling in unsaturated soil with incomplete physical constraints}.
\newblock \bibinfo{journal}{Computers and Geotechnics} \bibinfo{volume}{180}, \bibinfo{pages}{107091}.
\bibitem[{Njifenjou(2025)}]{njifenjou2025discrete}
\bibinfo{author}{Njifenjou, A.}, \bibinfo{year}{2025}.
\newblock \bibinfo{title}{Discrete maximum principle honored by conventional finite volume schemes for diffusion-convection-reaction problems: Proof with geometrical arguments}.
\newblock \bibinfo{journal}{London Journal of Research In Science: Natural and Formal} \bibinfo{volume}{25}, \bibinfo{pages}{43--57}.
\bibitem[{Or et~al.(2015)Or, Lehmann and Assouline}]{or2015}
\bibinfo{author}{Or, D.}, \bibinfo{author}{Lehmann, P.}, \bibinfo{author}{Assouline, S.}, \bibinfo{year}{2015}.
\newblock \bibinfo{title}{Natural length scales define the range of applicability of the {R}ichards equation for capillary flows}.
\newblock \bibinfo{journal}{Water Resources Research} \bibinfo{volume}{51}, \bibinfo{pages}{7130--7144}.
\bibitem[{Orouskhani et~al.(2023)Orouskhani, Sahoo, Agyeman, Bo and Liu}]{orouskhani_impact_2023}
\bibinfo{author}{Orouskhani, E.}, \bibinfo{author}{Sahoo, S.}, \bibinfo{author}{Agyeman, B.}, \bibinfo{author}{Bo, S.}, \bibinfo{author}{Liu, J.}, \bibinfo{year}{2023}.
\newblock \bibinfo{title}{Impact of sensor placement in soil water estimation: a real-case study}.
\newblock \bibinfo{journal}{Irrigation Science} \bibinfo{volume}{41}, \bibinfo{pages}{395--411}.
\bibitem[{Pichi et~al.(2024)Pichi, Moya and Hesthaven}]{pichi2024graph}
\bibinfo{author}{Pichi, F.}, \bibinfo{author}{Moya, B.}, \bibinfo{author}{Hesthaven, J.S.}, \bibinfo{year}{2024}.
\newblock \bibinfo{title}{A graph convolutional autoencoder approach to model order reduction for parametrized {PDEs}}.
\newblock \bibinfo{journal}{Journal of Computational Physics} \bibinfo{volume}{501}, \bibinfo{pages}{112762}.
\bibitem[{Pinkus(1999)}]{pinkus}
\bibinfo{author}{Pinkus, A.}, \bibinfo{year}{1999}.
\newblock \bibinfo{title}{Approximation theory of the {MLP} model in neural networks}.
\newblock \bibinfo{journal}{Acta Numerica} \bibinfo{volume}{8}, \bibinfo{pages}{143–195}.
\newblock \DOIprefix\doi{10.1017/S0962492900002919}.
\bibitem[{Raissi et~al.(2017)Raissi, Perdikaris and Karniadakis}]{raissi2017physics}
\bibinfo{author}{Raissi, M.}, \bibinfo{author}{Perdikaris, P.}, \bibinfo{author}{Karniadakis, G.E.}, \bibinfo{year}{2017}.
\newblock \bibinfo{title}{Physics informed deep learning (part i): Data-driven solutions of nonlinear partial differential equations}.
\newblock \bibinfo{journal}{arXiv preprint arXiv:1711.10561} .
\bibitem[{Raissi et~al.(2019)Raissi, Perdikaris and Karniadakis}]{raissi2019physics}
\bibinfo{author}{Raissi, M.}, \bibinfo{author}{Perdikaris, P.}, \bibinfo{author}{Karniadakis, G.E.}, \bibinfo{year}{2019}.
\newblock \bibinfo{title}{Physics-informed neural networks: A deep learning framework for solving forward and inverse problems involving nonlinear partial differential equations}.
\newblock \bibinfo{journal}{Journal of Computational Physics} \bibinfo{volume}{378}, \bibinfo{pages}{686--707}.
\bibitem[{Ranade et~al.(2021)Ranade, Hill and Pathak}]{ranade2021discretizationnet}
\bibinfo{author}{Ranade, R.}, \bibinfo{author}{Hill, C.}, \bibinfo{author}{Pathak, J.}, \bibinfo{year}{2021}.
\newblock \bibinfo{title}{Discretizationnet: A machine-learning based solver for {Navier--Stokes} equations using finite volume discretization}.
\newblock \bibinfo{journal}{Computer Methods in Applied Mechanics and Engineering} \bibinfo{volume}{378}, \bibinfo{pages}{113722}.
\bibitem[{Rathfelder and Abriola(1994)}]{massconservative}
\bibinfo{author}{Rathfelder, K.}, \bibinfo{author}{Abriola, L.M.}, \bibinfo{year}{1994}.
\newblock \bibinfo{title}{Mass conservative numerical solutions of the head-based {R}ichards equation}.
\newblock \bibinfo{journal}{Water Resources Research} \bibinfo{volume}{30}, \bibinfo{pages}{2579--2586}.
\bibitem[{Richards(1931)}]{richards}
\bibinfo{author}{Richards, L.A.}, \bibinfo{year}{1931}.
\newblock \bibinfo{title}{Capillary conduction of liquids through porous mediums}.
\newblock \bibinfo{journal}{Physics} \bibinfo{volume}{1}, \bibinfo{pages}{318--333}.
\bibitem[{Roth(2008)}]{roth2008}
\bibinfo{author}{Roth, K.}, \bibinfo{year}{2008}.
\newblock \bibinfo{title}{Scaling of water flow through porous media and soils}.
\newblock \bibinfo{journal}{European Journal of Soil Science} \bibinfo{volume}{59}, \bibinfo{pages}{125--130}.
\bibitem[{Saeedmonir et~al.(2024)Saeedmonir, Adeli and Khoei}]{saeedmonir2024multiscale}
\bibinfo{author}{Saeedmonir, S.}, \bibinfo{author}{Adeli, M.}, \bibinfo{author}{Khoei, A.}, \bibinfo{year}{2024}.
\newblock \bibinfo{title}{A multiscale approach in modeling of chemically reactive porous media}.
\newblock \bibinfo{journal}{Computers and Geotechnics} \bibinfo{volume}{165}, \bibinfo{pages}{105818}.
\bibitem[{{\v{S}}im{\r u}nek et~al.(2016){\v{S}}im{\r u}nek, Van~Genuchten and {\v{S}}ejna}]{vsimuunek2016recent}
\bibinfo{author}{{\v{S}}im{\r u}nek, J.}, \bibinfo{author}{Van~Genuchten, M.T.}, \bibinfo{author}{{\v{S}}ejna, M.}, \bibinfo{year}{2016}.
\newblock \bibinfo{title}{Recent developments and applications of the {HYDRUS} computer software packages}.
\newblock \bibinfo{journal}{Vadose Zone Journal} \bibinfo{volume}{15}.
\bibitem[{Sizikov et~al.(2011)}]{sizikov2011well}
\bibinfo{author}{Sizikov, V.S.}, et~al., \bibinfo{year}{2011}.
\newblock \bibinfo{title}{Well-posed, ill-posed, and intermediate problems with applications}.
\newblock \bibinfo{publisher}{De Gruyter}.
\bibitem[{Smith et~al.(2002)Smith, Smettem, Broadbridge and Woolhiser}]{smith_infiltration_2002}
\bibinfo{author}{Smith, R.E.}, \bibinfo{author}{Smettem, K.}, \bibinfo{author}{Broadbridge, P.}, \bibinfo{author}{Woolhiser, D.}, \bibinfo{year}{2002}.
\newblock \bibinfo{title}{Infiltration Theory for Hydrologic Applications}.
\newblock \bibinfo{publisher}{American Geophysical Union}.
\bibitem[{Song and Jiang(2023a)}]{songescape}
\bibinfo{author}{Song, Z.}, \bibinfo{author}{Jiang, Z.}, \bibinfo{year}{2023}a.
\newblock \bibinfo{title}{A data-driven modeling approach for water flow dynamics in soil}, in: \bibinfo{editor}{Kokossis, A.C.}, \bibinfo{editor}{Georgiadis, M.C.}, \bibinfo{editor}{Pistikopoulos, E.} (Eds.), \bibinfo{booktitle}{33rd European Symposium on Computer Aided Process Engineering}. \bibinfo{publisher}{Elsevier}. volume~\bibinfo{volume}{52} of \textit{\bibinfo{series}{Computer Aided Chemical Engineering}}, pp. \bibinfo{pages}{819--824}.
\bibitem[{Song and Jiang(2023b)}]{songdata}
\bibinfo{author}{Song, Z.}, \bibinfo{author}{Jiang, Z.}, \bibinfo{year}{2023}b.
\newblock \bibinfo{title}{A data-driven random walk approach for solving water flow dynamics in soil systems}, in: \bibinfo{booktitle}{Proceedings of Foundations of Computer-Aided Process Operations and Chemical Process Control (FOCAPO/CPC)}, pp. \bibinfo{pages}{1--6}.
\bibitem[{Song and Jiang(2025a)}]{escape_25_inverse}
\bibinfo{author}{Song, Z.}, \bibinfo{author}{Jiang, Z.}, \bibinfo{year}{2025}a.
\newblock \bibinfo{title}{A novel bayesian framework for inverse problems in precision agriculture}.
\newblock \bibinfo{journal}{Systems \& Control Transactions} \bibinfo{volume}{4}, \bibinfo{pages}{246--251}.
\bibitem[{Song and Jiang(2025b)}]{escape_25_fractional}
\bibinfo{author}{Song, Z.}, \bibinfo{author}{Jiang, Z.}, \bibinfo{year}{2025}b.
\newblock \bibinfo{title}{A physics-based, data-driven numerical framework for anomalous diffusion of water in soil}.
\newblock \bibinfo{journal}{Systems \& Control Transactions} \bibinfo{volume}{4}, \bibinfo{pages}{2391--2397}.
\bibitem[{Su et~al.(2022)Su, Zhang, Zou, Zhao, Zhang and Su}]{su2022numerical}
\bibinfo{author}{Su, X.}, \bibinfo{author}{Zhang, M.}, \bibinfo{author}{Zou, D.}, \bibinfo{author}{Zhao, Y.}, \bibinfo{author}{Zhang, J.}, \bibinfo{author}{Su, H.}, \bibinfo{year}{2022}.
\newblock \bibinfo{title}{Numerical scheme for solving the {Richard’s} equation based on finite volume model with unstructured mesh and implicit dual-time stepping}.
\newblock \bibinfo{journal}{Computers and Geotechnics} \bibinfo{volume}{147}, \bibinfo{pages}{104768}.
\bibitem[{Sun et~al.(2024)Sun, Li, Jiao, Liu, Yang, Lyu and Lin}]{sun2024simulation}
\bibinfo{author}{Sun, Y.}, \bibinfo{author}{Li, D.}, \bibinfo{author}{Jiao, L.}, \bibinfo{author}{Liu, Z.}, \bibinfo{author}{Yang, Y.}, \bibinfo{author}{Lyu, J.}, \bibinfo{author}{Lin, S.}, \bibinfo{year}{2024}.
\newblock \bibinfo{title}{Simulation of saturated--unsaturated seepage problems via the virtual element method}.
\newblock \bibinfo{journal}{Computers and Geotechnics} \bibinfo{volume}{171}, \bibinfo{pages}{106326}.
\bibitem[{Tracy(2006)}]{tracy2006clean}
\bibinfo{author}{Tracy, F.T.}, \bibinfo{year}{2006}.
\newblock \bibinfo{title}{Clean two-and three-dimensional analytical solutions of {Richards'} equation for testing numerical solvers}.
\newblock \bibinfo{journal}{Water Resources Research} \bibinfo{volume}{42}.
\bibitem[{Van~Genuchten(1980)}]{vanGenuchten1980}
\bibinfo{author}{Van~Genuchten, M.T.}, \bibinfo{year}{1980}.
\newblock \bibinfo{title}{A closed-form equation for predicting the hydraulic conductivity of unsaturated soils}.
\newblock \bibinfo{journal}{Soil Science Society of America Journal} \bibinfo{volume}{44}, \bibinfo{pages}{892--898}.
\bibitem[{Venkatesan et~al.(2021)Venkatesan, Droniou, Roy, Robert and Zhou}]{venkatesan2021application}
\bibinfo{author}{Venkatesan, S.}, \bibinfo{author}{Droniou, J.}, \bibinfo{author}{Roy, I.}, \bibinfo{author}{Robert, D.}, \bibinfo{author}{Zhou, A.}, \bibinfo{year}{2021}.
\newblock \bibinfo{title}{Application of diffusion-advection equations to in-field monitoring of soil suction profiles}.
\newblock \bibinfo{journal}{Computers and Geotechnics} \bibinfo{volume}{139}, \bibinfo{pages}{104329}.
\bibitem[{Vogel and Ippisch(2008)}]{Vogel2008}
\bibinfo{author}{Vogel, H.J.}, \bibinfo{author}{Ippisch, O.}, \bibinfo{year}{2008}.
\newblock \bibinfo{title}{Estimation of a critical spatial discretization limit for solving {R}ichards' equation at large scales}.
\newblock \bibinfo{journal}{Vadose Zone Journal} \bibinfo{volume}{7}, \bibinfo{pages}{112--114}.
\bibitem[{Zarba(1988)}]{zarba1988thesis}
\bibinfo{author}{Zarba, R.}, \bibinfo{year}{1988}.
\newblock \bibinfo{title}{A Numerical Investigation of Unsaturated Flow}.
\newblock \bibinfo{publisher}{Massachusetts Institute of Technology, Department of Civil Engineering}, \bibinfo{address}{Cambridge, MA}.
\bibitem[{Zeidler(1986)}]{zeidler2013nonlinear}
\bibinfo{author}{Zeidler, E.}, \bibinfo{year}{1986}.
\newblock \bibinfo{title}{Nonlinear functional analysis and its applications. I. Fixed-Point Theorems}.
\newblock \bibinfo{publisher}{Springer Science \& Business Media}.
\bibitem[{Zhang et~al.(2018)Zhang, Zhu, Shi and Fatahi}]{zhang2018long}
\bibinfo{author}{Zhang, C.C.}, \bibinfo{author}{Zhu, H.H.}, \bibinfo{author}{Shi, B.}, \bibinfo{author}{Fatahi, B.}, \bibinfo{year}{2018}.
\newblock \bibinfo{title}{A long term evaluation of circular mat foundations on clay deposits using fractional derivatives}.
\newblock \bibinfo{journal}{Computers and Geotechnics} \bibinfo{volume}{94}, \bibinfo{pages}{72--82}.
\bibitem[{Zheng et~al.(2025)Zheng, Li, Wan, Yang and Zhu}]{zheng2025wave}
\bibinfo{author}{Zheng, S.}, \bibinfo{author}{Li, W.}, \bibinfo{author}{Wan, Y.}, \bibinfo{author}{Yang, Z.}, \bibinfo{author}{Zhu, S.}, \bibinfo{year}{2025}.
\newblock \bibinfo{title}{Wave propagation in an ocean site considering fractional viscoelastic constitution of porous seabed}.
\newblock \bibinfo{journal}{Computers and Geotechnics} \bibinfo{volume}{180}, \bibinfo{pages}{107098}.
\bibitem[{Zhou et~al.(2013)Zhou, Liu and He}]{zhou2013richards}
\bibinfo{author}{Zhou, J.}, \bibinfo{author}{Liu, F.}, \bibinfo{author}{He, J.H.}, \bibinfo{year}{2013}.
\newblock \bibinfo{title}{On {Richards’} equation for water transport in unsaturated soils and porous fabrics}.
\newblock \bibinfo{journal}{Computers and Geotechnics} \bibinfo{volume}{54}, \bibinfo{pages}{69--71}.
\bibitem[{Zhu et~al.(2019)Zhu, Wu, Shen and Huang}]{zhu2019improved}
\bibinfo{author}{Zhu, S.}, \bibinfo{author}{Wu, L.}, \bibinfo{author}{Shen, Z.}, \bibinfo{author}{Huang, R.}, \bibinfo{year}{2019}.
\newblock \bibinfo{title}{An improved iteration method for the numerical solution of groundwater flow in unsaturated soils}.
\newblock \bibinfo{journal}{Computers and Geotechnics} \bibinfo{volume}{114}, \bibinfo{pages}{103113}.

\end{thebibliography}
\end{document}